\documentclass[12pt,twoside,reqno,openany]{amsart}

\usepackage{amsbsy,amscd,amsfonts,amsmath,amssymb,amsthm,color,
fancybox,fancyhdr,footmisc,graphics,graphicx,ifthen,mathrsfs,
multicol,pdfpages,rotating,times,wasysym}
\usepackage[dvipsnames,svgnames,x11names]{xcolor}

\usepackage[all]{xy}
\usepackage[utf8]{inputenc}
\usepackage[T1]{fontenc}
\usepackage{mathtools}
\newtagform{EngelLie}[\scriptstyle]{$}{$}
\makeatletter\newcommand{\leqnomode}{\tagsleft@true}
\newcommand{\reqnomode}{\tagsleft@false}\makeatother

\newtheorem{Theorem}[equation]{Theorem}

\newtheorem{Proposition}[equation]{Proposition}

\newtheorem{Lemma}[equation]{Lemma}

\newtheorem{Corollary}[equation]{Corollary}

\theoremstyle{definition}

\newtheorem{Definition}[equation]{D\'efinition}

\newtheorem{Example}[equation]{Example}
\newtheorem{Remark}[equation]{Remark}

\newtheorem{Problem}[equation]{Problem}

\newcommand{\C}{\mathbb{C}}

\renewcommand{\H}{\mathbb{H}}

\renewcommand{\P}{\mathbb{P}}

\newcommand{\R}{\mathbb{R}}

\newcommand{\KK}{\text{\sc k}}

\definecolor{blue}{cmyk}{1.,1.,0.,0.63}
\definecolor{red}{cmyk}{0.,1.,1.,0.63}
\definecolor{green}{cmyk}{1.,0.,1.,0.63}
\definecolor{black}{cmyk}{1.,1.,1.,1.}

\newcommand{\blue}{\textcolor{blue}}
\newcommand{\green}{\textcolor{green}}

\makeatletter
\renewcommand{\@fnsymbol}[1]
{\ensuremath{\ifcase#1\or $*$\or $**$\or $***$\or $****$\or $*****$
\else\@ctrerr\fi}}
\makeatother

\numberwithin{equation}{section}

\newcommand{\Section}[1]{
\renewcommand{\thesection}{\bf\arabic{section}}
\section{#1}
\renewcommand{\thesection}{\arabic{section}}}

\newcommand{\style}[1]{\text{\footnotesize{\sf #1}}}

\newcommand{\stylesmall}[1]{{\sf #1}}

\renewcommand{\arccos}{\style{arccos}}

\newcommand{\arccosh}{\style{arccosh}}

\renewcommand{\arcsin}{\style{arcsin}}

\newcommand{\arcsinh}{\style{arcsinh}}

\newcommand{\Aut}{\style{Aut}}

\renewcommand{\cos}{\style{cos}}
\newcommand{\cossmall}{\stylesmall{cos}}

\renewcommand{\cosh}{\style{cosh}}

\renewcommand{\det}{\style{det}}

\newcommand{\dist}{\style{dist}}

\newcommand{\Ellsmall}{\stylesmall{Ell}}

\newcommand{\End}{\style{End}}

\newcommand{\Euclsmall}{\stylesmall{Eucl}}

\newcommand{\Flatsmall}{\stylesmall{Flat}}

\newcommand{\Hypsmall}{\stylesmall{Hyp}}

\newcommand{\Idsmall}{\stylesmall{Id}}

\renewcommand{\Im}{\style{Im}}
\newcommand{\Imsmall}{\stylesmall{Im}}

\renewcommand{\lim}{\style{lim}}

\renewcommand{\log}{\style{log}}

\renewcommand{\max}{\style{max}}

\renewcommand{\mod}{\style{mod}}

\renewcommand{\Re}{\style{Re}}
\newcommand{\Resmall}{\stylesmall{Re}}

\renewcommand{\sin}{\style{sin}}
\newcommand{\sinsmall}{\stylesmall{sin}}

\renewcommand{\sinh}{\style{sinh}}
\newcommand{\sinhsmall}{\stylesmall{sinh}}

\newcommand{\isqrt}{{\scriptstyle{\sqrt{-1}}}}

\newcommand{\medcup}{\mathbin{\scalebox{1.5}{\ensuremath{\cup}}}}

\newcommand{\smallbullet}{{\scriptscriptstyle{\bullet}}}

\newcommand{\vf}{\vfill\end{document}}

\begin{document}

\title[Non-umbilicity of hyperbolic surfaces Grauert tubes]{
Nonvanishing of Cartan CR curvature
\\
on boundaries of Grauert tubes 
\\
around hyperbolic surfaces}

\author{Wei Guo Foo}

\address{Hua Loo-Keng center for mathematical sciences, Academy of Mathematics and Systems Sciences, Chinese Academy of Sciences, Beijing, China}

\email{fooweiguo@hotmail.com}

\author{Jo\"{e}l Merker}

\address{Laboratoire de Math\'{e}matiques d'Orsay, 
B\^{a}timent 307, Universit\'{e} Paris-Sud, CNRS, Universit\'{e} Paris-Saclay, 91405 Orsay Cedex, France}

\email{joel.merker@math.u-psud.fr, the-anh.ta@u-psud.fr}

\author{The-Anh Ta}

\subjclass[2000]{32V25, 32V40}

\date{\today}

\keywords{Poincar\'e metric, Cartan CR-curvature, 
Grauert tube, Umbilical points}

\begin{abstract}
We show that the boundaries of thin strongly pseudoconvex Grauert
tubes, with respect to the Guillemin-Stenzel K\"{a}hler metric
canonically associated with the Poincar\'e metric on closed hyperbolic
real-analytic surfaces, has nowhere vanishing Cartan
CR-curvature. This result provides a wealth of examples of compact
$3$-dimensional Levi nondegenerate CR manifolds having no CR-umbilical
point.

We provide two proofs utilizing two recent formulas for determining the
Cartan CR-curvature of any local $\mathcal{C}^6$-smooth hypersurfaces
in $\C^2$. One was obtained in 2012 by the second
named author joint with Sabzevari, and it is an {\em expanded}
explicit formula, valid for locally graphed hypersurfaces, containing
millions of terms. The other formula, which we published in 2018 when
studying Webster's ellipsoidal hypersurfaces, is not expanded, but
more suitable for calculations with a hypersurface in $\mathbb{C}^2$
that is represented as the zero locus of some {\em
implicit}\,\,---\,\,but `simple' in some sense, {\em e.g.}
quadratic\,\,---\,\,defining function.

We also discuss Grauert tubes constructed with respect to extrinsic
metrics depending
on embeddings in complex surfaces, 
together with a certain combinatorics of product metrics.
\end{abstract}

\maketitle

\Section{\bf Introduction}
\label{introduction}

The equivalence problem for local real-analytic hypersurfaces with
respect to local biholomorphisms in $\mathbb{C}^2$ was first studied
by Poincar\'{e} \cite{Poincare-1907}, and was later solved by
Cartan~\cite{Cartan-1932} with the introduction of the so-called
method of equivalence.  The theory was later developed in 
$\mathbb{C}^{n+1}$ by Chern and
Moser~\cite{Chern-Moser-1974}, and resulted in the set up of invariant
CR-curvatures, called Cartan curvatures in complex dimension 2, and
Hachtroudi-Chern curvatures when $n \geqslant 3$. 

For a long time,
little was known about these curvatures due to their high
computational complexity. Nonetheless, Webster~\cite{Webster-2000},
and later Huang and Ji~\cite{Huang-Ji-2007} were able to investigate
the case of real ellipsoidal hypersurfaces. In recent years, new
variants and explicit formulas 
({\em see}~\cite{Ebenfelt-Duong-Zaitsev-2016,
Merker-Sabzevari-2012, Merker-Sabzevari-2014,
Foo-Merker-Ta-crmath})
made it possible to determine the vanishing
locus of the Cartan curvatures for new classes of $3$-dimensional CR
manifolds. For instance, we were able to find a whole
explicit curve of points of vanishing Cartan curvature on general
ellipsoids in $\mathbb{C}^2$ in~\cite{Foo-Merker-Ta-crmath}.

In their landmark paper~\cite[p.~247]{Chern-Moser-1974}, 
Chern and Moser raised the following

\begin{Problem}
{\sl Are there compact strictly pseudoconvex 
hypersurfaces $M^3 \subset \C^2$ without
CR-umbilical points?  Are there such manifolds diffeomorphic to the
sphere $S^3 \subset \C^2$.}
\end{Problem}

It is well known that a standard $2$-torus in $\R^3$ has no
Riemannian-umbilic point.  Similarly, it is not difficult to verify
(\cite{Ebenfelt-Duong-Zaitsev-2016}) that the boundaries of thin
Grauert tubes around the flat 2-dimensional torus $\mathbb{T}^2 = S^1
\times S^1 \subset \mathbb{C}^2$ have empty CR-umbilical locus.
Thus, a topological restriction like $M^3 \cong S^3$ must be
assumed.

In this paper, we are interested in the question of whether a similar
phenomenon holds for higher genus surfaces. Let therefore $S$ be a
closed compact real-analytic ($\mathcal{C}^\omega$) surface of genus
$\geqslant 2$ which is hyperbolic in the sense that its universal
cover is the unit disc $\mathbb{D} \subset \mathbb{C}$. As a special
case of a theorem of Bruhat and Whitney~\cite{Bruhat-Whitney-1959} in
dimension 2, $S$ admits an extrinsic complexification, namely there
exists a complex manifold $M^c$ of complex dimension 2, together with
an analytic totally real embedding of $S$ into $M^c$. Moreover, the
work~\cite{Guillemin-Stenzel-1991} of Guillemin and Stenzel provides a
canonical K\"{a}hler potential $\rho$ defined in a small neighborhood
of $S$ in $M^c$ ({\em see} 
Section~{\ref{Khaler-potential-Grauert-tubes}} 
below). In particular, for each
$\varepsilon$ with $0 < \varepsilon \leqslant \varepsilon_0 \ll 1$,
the set $\Omega_\varepsilon := \rho^{-1} \big( [ 0, \varepsilon)
\big)$, called the Grauert tube of radius $\varepsilon$ around $S$,
has strongly pseudoconvex $\mathcal{C}^\omega$ boundary
$M_{\varepsilon} := \rho^{-1} (\varepsilon)$ contained in the
complex surface $M^c$, to which Cartan's method of equivalence
applies. Our main result is the following.

\begin{Theorem}
\label{main-theorem}
There exists $0 < \varepsilon_0 \ll 1$ such that for every
$\varepsilon$ with $0 < \varepsilon \leqslant \varepsilon_0$, 
the real and imaginary parts 
of the primary complex
Cartan curvature vanish nowhere on the boundary of
$M_\varepsilon$.
\end{Theorem}

Equivalently: 

\begin{Corollary}
The boundaries of these $M_\varepsilon$ have no CR-umbilical point.\qed
\end{Corollary}

So far, our construction of the Grauert tubes $\Omega_\varepsilon$
take a complete intrisic point of view, since the Guillemin-Stenzel
potential is obtained only from a given intrinsic metric on the
surface $S$. It is then natural to look at the Grauert tubes from an
extrinsic point of view, that is we consider the surface $S$ as being
totally really embedded in a given (local) complex surface equipped
with a given metric. Already in the case of a torus embedded in the
standard $\mathbb{C}^2$, the extrinsic contruction will provide
several new examples of compact hypersurfaces without CR-umbilical
points ({\em see} Example~{\ref{Example-flat-products}}). Further
constructions in this vein are provided in
Section~{\ref{Grauert-extrinsic-metrics}}.

\smallskip

This paper is organized as follows. In
Section~{\ref{Khaler-potential-Grauert-tubes}}, we recall the
construction of the canonical K\"{a}hler potential of Guillemin and
Stenzel in \cite{Guillemin-Stenzel-1991}, and we find an explicit
formula for the potential in the case of hyperbolic surfaces.
Section~{\ref{round-sphere-flat-torus}} discusses two standard
examples of complexification of the round sphere and the flat
torus. In Section~{\ref{semi-global-tube}}, we work out the defining
function for the Grauert tube around the Poincar\'{e} upper
half-plane. The formula then will be used in
Section~{\ref{calculation-complex-curvature}} to calculate the Cartan
curvatures on the boundaries of Grauert tubes of hyperbolic surfaces
by explicit expressions given in \cite{Merker-Sabzevari-2012} and
\cite{Foo-Merker-Ta-crmath}, and to show that the Cartan curvatures do
not vanish for small enough radii.  Section~{\ref{transfer-genus-g}}
explains in details how nonvanishing of the Cartan curvature on the
boundary of Grauert tubes around hyperbolic surfaces can be deduced
from the calculations in
Section~{\ref{calculation-complex-curvature}}. Finally, in
Section~{\ref{Grauert-extrinsic-metrics}}, we discuss some extrinsic
constructions of Grauert tubes based on product metrics.

\Section{\bf The Canonical K\"{a}hler Potential on Grauert Tubes}
\label{Khaler-potential-Grauert-tubes}

For any compact real-analytic ($\mathcal{C}^\omega$) manifold $M$ of
dimension $n \geqslant 1$, Bruhat and Whitney showed
in~\cite{Bruhat-Whitney-1959} that there exists an $n$-dimensional
complex manifold $M^c$, and a real-analytic embedding $M
\hookrightarrow M^c$ which is totally real, {\em i.e.} such that the
real tangent spaces to $M$ contain no complex lines in the complex
tangent spaces to $M^c$.  The $\mathcal{C}^\omega$ changes of charts
$\R^n \ni x \longmapsto x' = \varphi(x) \in \R^n$ for $M$, where $x =
(x_1, \dots, x_n)$, become $\C^n \ni z \longmapsto z' = \varphi(z) \in
\C^n$, where $z = x + \isqrt\, y \in \C^n$, and where $\varphi(z)$
means substituting $z$ for $x$ in the punctual convergent power series
of $\varphi$, giving the complex manifold structure of $M^c$.  The
Taylor coefficients of such $\mathcal{C}^\omega$ diffeomorphisms
$\varphi = \overline{\varphi}$ are real, the complex
conjugation $z \longmapsto \overline{z}$ transfers coherently as
$\overline{z}' = \varphi(\overline{z})$, which shows that $M = {\rm
Fix}(\sigma)$ is the set of fixed points of the antiholomorphic involution $\sigma \colon M^c \longrightarrow M^c$ obtained from $z \longmapsto \overline{z}$ in any chart.

Also by substituting $z$ for $x$ in power series, every
$\mathcal{C}^\omega$ function $f \colon M \longrightarrow \R$ extends
uniquely as a holomorphic function $f^c \colon U^c \longrightarrow \C$
with $f^c \big\vert_M = f$, in some open neighborhood $U^c$ of $M$ in $M^c: M \subset U^c
\subset M^c$, and $f \big\vert_M \equiv 0$ if and only if $f^c \equiv
0$ in some subneighborhood $V^c:$ $M \subset V^c \subset U^c$.

According to Grauert~\cite{Grauert-1958}, there exists a
$\mathcal{C}^\infty$ strictly plurisubharmonic function $\rho \colon
U^c \longrightarrow [0,1)$ defined in some open neighborhood $U^c$ of $M$ in $M^c,$ with $ \rho\circ \sigma = \rho, M =
\rho^{-1}(0), d\rho \big\vert_M \equiv 0$, and such that
$\rho$ has no critical point in $V^c \backslash M$, for
some subneighborhood $V^c: M \subset V^c \subset U^c$.
Hence for all small enough $\varepsilon: 0 < \varepsilon \leqslant 
\varepsilon_0 \ll 1$,
the domain $\Omega_\varepsilon = \{ \rho < \varepsilon\}$, a tubular
neighborhood of $M$ in $M^c$, has $\mathcal{C}^\infty$
strictly pseudoconvex boundary $M_\varepsilon = \{\rho =
\varepsilon\}$, and is called the {\sl Grauert tube} of radius
$\varepsilon$ around $M$.

When the manifold $M$ is equipped with some $\mathcal{C}^\omega$
Riemannian metric $g$, Guillemin and Stenzel gave in
\cite{Guillemin-Stenzel-1991} a very elegant construction of such a
strictly plurisubharmonic function
\[
\rho 
\,=\, 
\rho_g
\colon
\ \ \
M^c 
\longrightarrow
[0,1)
\]
uniquely associated to $g$ that will be called the {\sl canonical
K\"{a}hler potential} on $M^c$. Their construction can be
summarized as follows.

Embed $M \hookrightarrow M \times M$ by $x \longmapsto (x,x)$ and let
$W$ be an open neighborhood of $M$ in $M \times M$. If $W$ is thin
enough, for any pair $(x,u) \in W$, the local uniqueness and distance
minimizing properties of geodesics with respect to $g$ guarantees that
$\dist_g(x,u)$ is the $g$-length of the geodesic from $x$ to $u$, and
an inspection of the $g$-length formula convinces that the (symmetric)
{\em squared} distance function:
\[
f(x,u)
\,:=\,
\big(
\dist_g(x,u)
\big)^2
\eqno
{\scriptstyle{(x,\,u\,\in\,W)}}
\]
is $\mathcal{C}^\omega$, hence can be complexified.

Since in local coordinates, we will denote $x
= (x_1, \dots, x_n)$ and $u = (u_1, \dots, u_n)$ in $\R^n$ and
introduce $z := x + \isqrt\, y$ with $w := u + \isqrt\, v$ in $\C^n$,
let us denote a pair of points in the global abstract product
similarly as $(z,w) \in M^c \times M^c$, and let us abbreviate $\sigma
\colon M^c \longrightarrow M^c$ as $z \longmapsto \overline{z}$.
Also, let us use the embedding:
\[
M^c
\,\ni\,
z
\,\longmapsto\,
\big(z,\overline{z}\big)
\,\in\,M^c\times M^c,
\]
compatible with $x \longmapsto (x,x)$ which makes $M^c$ totally real
in $M^c \times M^c$, and let $W^c$ be a thin open neighborhood of
$M^c$ in $M^c \times M^c$ invariant under the conjugation $(z, w) \longmapsto \big(
\overline{w}, \overline{z} \big)$ and satisfying $W = W^c \cap (M \times
M)$.

\begin{center}
\begin{picture}(0,0)%
\includegraphics{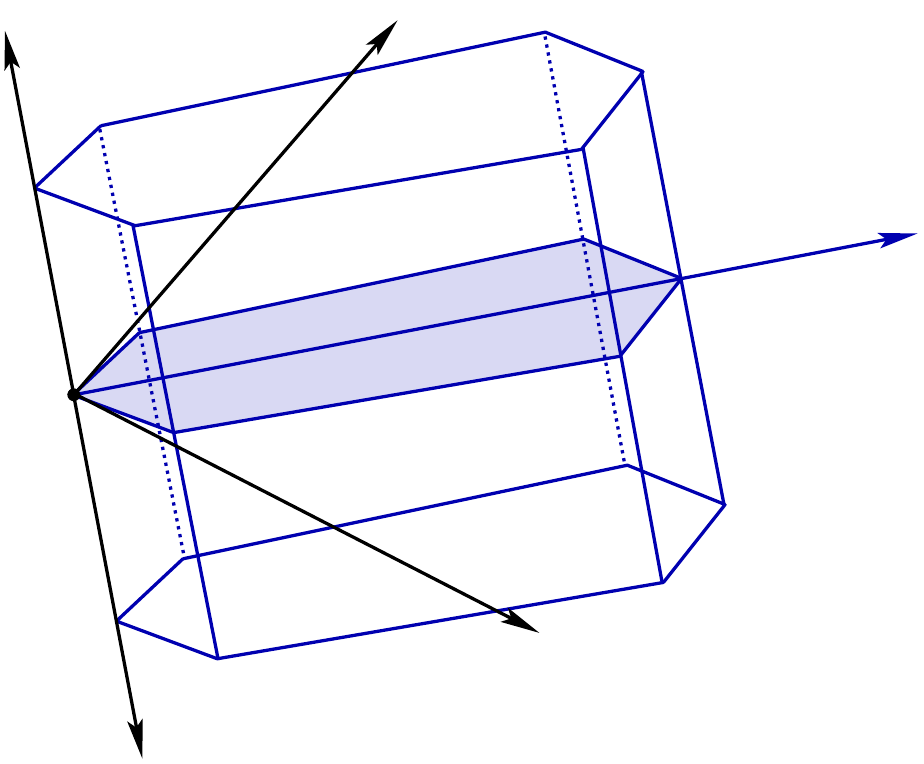}%
\end{picture}%
\setlength{\unitlength}{4144sp}%
\begingroup\makeatletter\ifx\SetFigFont\undefined%
\gdef\SetFigFont#1#2#3#4#5{%
  \reset@font\fontsize{#1}{#2pt}%
  \fontfamily{#3}\fontseries{#4}\fontshape{#5}%
  \selectfont}%
\fi\endgroup%
\begin{picture}(4218,3477)(1554,-3188)
\put(4024,-2580){\makebox(0,0)[lb]{\smash{{\SetFigFont{9}{10.8}{\familydefault}{\mddefault}{\updefault}{\color[rgb]{0,0,0}$M$}%
}}}}
\put(2224,-3116){\makebox(0,0)[lb]{\smash{{\SetFigFont{9}{10.8}{\familydefault}{\mddefault}{\updefault}{\color[rgb]{0,0,0}$M^c$}%
}}}}
\put(4039,-2310){\makebox(0,0)[lb]{\smash{{\SetFigFont{9}{10.8}{\familydefault}{\mddefault}{\updefault}{\color[rgb]{0,0,0}\blue{$W^c$}}%
}}}}
\put(3618,-101){\makebox(0,0)[lb]{\smash{{\SetFigFont{9}{10.8}{\familydefault}{\mddefault}{\updefault}{\color[rgb]{0,0,0}\blue{$W^c$}}%
}}}}
\put(1627, 83){\makebox(0,0)[lb]{\smash{{\SetFigFont{9}{10.8}{\familydefault}{\mddefault}{\updefault}{\color[rgb]{0,0,0}$M^c$}%
}}}}
\put(3196,-1154){\makebox(0,0)[lb]{\smash{{\SetFigFont{9}{10.8}{\familydefault}{\mddefault}{\updefault}{\color[rgb]{0,0,0}\blue{$W$}}%
}}}}
\put(3241,-1409){\makebox(0,0)[lb]{\smash{{\SetFigFont{9}{10.8}{\familydefault}{\mddefault}{\updefault}{\color[rgb]{0,0,0}\blue{$W$}}%
}}}}
\put(3392,189){\makebox(0,0)[lb]{\smash{{\SetFigFont{9}{10.8}{\familydefault}{\mddefault}{\updefault}{\color[rgb]{0,0,0}$M$}%
}}}}
\end{picture}%

\end{center}

Then $f(x,u)$ complexifies as $f^c(z,w)$ defined and holomorphic for
$(z,w) \in W^c$, with $f^c \big\vert_M \equiv f$ and enjoys the
symmetry $f^c(w,z) = f^c(z,w)$. Furthermore, the reality condition
$\overline{f(x,u)} = f(x,u)$ of $f$ yields via complexification:
\[
\overline{f^c(z,w)}
\,\equiv\,
f^c\big(\overline{z},\overline{w}\big),
\]
hence putting $w := \overline{z}$, and using the symmetry,
we see the reality:
\[
\overline{f^c(z,\overline{z})}
\,\equiv\,
f^c\big(\overline{z},z\big)
\,\equiv\,
f^c(z,\overline{z}).
\]

\begin{Proposition}
{\rm (\cite{Guillemin-Stenzel-1991}, p.~565)}
The real-valued function $f^c\big( z, \overline{z} \big)$ is
equal to $0$ on $M \hookrightarrow M^c \times M^c$ and  
takes values $<0$ outside $M$.
\end{Proposition}

So in $W^c \big\backslash \{ f^c = 0\}$, the square root $\sqrt{f^c}$
is $2:1$-valued, and the canonical K\"ahler potential $\rho = \rho_g$
is defined to be:
\leqnomode\usetagform{default}
\begin{align}
\label{rho-f-c}
\rho
\,:=\,
-f^c,
\end{align}
so that $\sqrt{\rho}$ is well defined in $\R_+$.

Finally, a consequence of Gauss' orthogonality lemma
(\cite{Guillemin-Stenzel-1991}, p.~564) 
which provides the annihilation:
\[
0
\,\equiv\,
\det\,
\bigg(
\frac{\partial^2\sqrt{f}}{\partial x_i\partial y_j}
(x,y)
\bigg)
\eqno
{\scriptstyle{(\forall\,(x,y)\,\in\,W\backslash M)}},
\]
yields via complexification the Monge-Amp\`ere equation:
\[
0
\,\equiv\,
\det\,
\bigg(
\frac{\partial^2\sqrt{\rho}}{\partial z_i\partial w_j}
(z,w)
\bigg)
\eqno
{\scriptstyle{(\forall\,(z,w)\,\in\,W^c\backslash M^c)}}.
\]
In \cite{Guillemin-Stenzel-1991}, Guillemin and Stenzel established the {\em uniqueness} of the
K\"ahler metric $\omega := \isqrt\, \partial \overline{\partial}
\rho_g$ on $M^c$ satisfying this and restricting to $g = \omega
\big\vert_M$ on $M$.

Of particular interest to us is the computational fact that $\rho =
\rho_g$ has explicit, workable expressions once $g$ is given,
especially in the case of surfaces.

\Section{\bf Two Examples: Round Sphere and Flat Torus}
\label{round-sphere-flat-torus}

\begin{Example}\cite[Section 4]{Guillemin-Stenzel-1991}
Consider $M : = \mathbb{S}^2$ to be the $2$-dimensional sphere:
\[
\mathbb{S}^2
\,:=\,
\big\{
(x_1,x_2,x_3) 
\in 
\mathbb{R}^3\colon\,
x_1^2+x_2^2+x_3^2=1
\big\},
\]
equipped with the standard round metric, whence the 
squared geodesic
distance between two points $x, y \in \mathbb{S}^2$ is:
\[
f
\big(
x,y
\big) 
\,=\,
\Big(
2\,\arcsin\,
\big(
{\textstyle{\frac{1}{2}}}
\sqrt{(x_1-y_1)^2+(x_2-y_2)^2+(x_3-y_3)^2}
\big)
\Big)^2.
\]
The Bruhat-Whitney complexification of $\mathbb{S}^2$ 
can be represented extrinsically as: 
\[
(\mathbb{S}^2)^c
\,:=\,
\big\{
(z_1,z_2,z_3)\in\mathbb{C}^3
\colon\,
z_1^2+z_2^2+z_3^2=1
\big\},
\]
and on it, we have the useful relation:
\[
(\Im\,z_1)^2+(\Im\,z_2)^2+(\Im\,z_3)^2
\,=\,
\big(z_1 \overline{z}_1 
+
z_2 \overline{z}_2
+
z_3 \overline{z}_3-1
\big)\big/2.
\]

The complexification of $f$ is:
\[
f^c
\big(
z,w
\big)
\,=\,
\Big(
2\,\arcsin\,
\big(
{\textstyle{\frac{1}{2}}}
\sqrt{
(z_1-w_1)^2+(z_2-w_2)^2+(z_3-w_3)^2}
\big)
\Big)^2,
\]
hence letting $w := \overline{z}$ and using the two 
identities:
\[
\arcsin\,\big(\isqrt\, t\big)
\,=\,
\isqrt\,\arcsinh (t),
\ \ \ \ \ \ \ \ \ \ \ \ \
\ \ \ \ \ \ \ \ \ \ \ \ \
2\,\arcsinh\,t 
\,=\,
\arccosh\,
\big(
1 + 2\,t^2
\big),
\]
we get:
\[
\aligned
f^c
\big(
z,\overline{z}
\big)
&
\,=\,
\bigg(
2\,\arcsin\,
\Big(
\pm\isqrt\,
\sqrt{
(\Im\,z_1)^2 +(\Im\,z_2)^2+(\Im\,z_3)^2}
\Big)
\bigg)^2
\\
&
\,=\,
\bigg(
\pm\,2\,\isqrt\,
\arcsinh\,
\Big(
{\textstyle{
\sqrt{
\big(z_1\overline{z}_1+z_2\overline{z}_2+z_3\overline{z}_3-1\big)
\big/2}}}
\Big)
\bigg)^2
\\
&
\,=\,
-\,
\Big(
\arccosh\,
\big(
z_1\overline{z}_1+z_2\overline{z}_2+z_3\overline{z}_3
\big)
\Big)^2,
\endaligned
\]
whence, coming back to the
definition~{\thetag{\ref{rho-f-c}}} of $\rho := -\, f^c$, we obtain:
\leqnomode\usetagform{default}
\begin{align}
\rho
\big(
z,\overline{z}
\big) 
\,=\,
\Big(
\arccosh\,
\big(
z_1\overline{z}_1+z_2\overline{z}_2+z_3\overline{z}_3
\big)
\Big)^2.
\end{align}
\end{Example}

\begin{Example}\label{flat-torus-example}
\cite[Section 3]{Ebenfelt-Duong-Zaitsev-2016}
Consider $M := \mathbb{T}^2 = \mathbb{R}^2 \big/ (2\pi \mathbb{Z}^2)$
to be the flat torus. Its complexification is $M^c := \mathbb{C}^2
\big/ (2\pi \mathbb{Z}^2)$. The geodesic distance between two close
points on $\mathbb{T}^2$ is computed along straight lines within the
flat universal cover $\big( \mathbb{R}^2,d_\Euclsmall \big)$. So, in a
fundamental domain for $\mathbb{T}^2$ on $\mathbb{R}^2$, the squared
distance and its complexification are
\[
\aligned
f\big(
(x_1,x_2),(y_1,y_2)
\big) 
&
\,=\,
(x_1-y_1)^2+(x_2-y_2)^2,
\\
f^c
\big(
(z_1,z_2),(w_1,w_2)
\big) 
&
\,=\,
(z_1-w_1)^2+(z_2-w_2)^2,
\endaligned
\]
hence letting $(w_1,w_2)=(\overline{z}_1,\overline{z}_2)^c$,
we get by the definition~\thetag{\ref{rho-f-c}} of $\rho := 
-\, f^c$:
\leqnomode\usetagform{default}
\begin{align}
\rho\big(
z,\overline{z}
\big)
\,=\,
4\,
(\Im\,z_1)^2
+
4\,
(\Im\,z_2)^2.
\end{align}
\end{Example}

\Section{\bf Semi-global Grauert Tube Around Poincaré's Upper 
Half-Plane}
\label{semi-global-tube}

For our purpose, we need to find the K\"ahler potential $\rho$ locally
on the Bruhat-Whitney complexification of any compact
$\mathcal{C}^\omega$ surface $S$ of genus $\geqslant 2$. When $S$ is
viewed as a Riemann surface, the uniformization
theorem (\cite[Chap.~27]{Forster-1991}) 
states that its universal cover is the
upper half-plane $\mathbb{H} = \{ z \in \C \colon\,
\Im(z) > 0\}$, and that:
\[
S
\,\cong\,
\mathbb{H}\big/\pi_1(S).
\]
We will then transfer geometric objects from $\H$ to $S$.

But in this section, our calculations will be done entirely in $\H =
\{ (x_1,x_2) \in \R\colon\, x_2 > 0\}$, viewed as a {\em real}
$\mathcal{C}^\omega$ surface equipped with the Poincaré metric $ds^2 =
\frac{dx_1^2 + dx_2^2}{x_2^2}$.  Since the 
squared Poincar\'e distance between
two points $(x_1, x_2)$ and $(y_1, y_2)$ of $\mathbb{H}$, with $x_2,
y_2 > 0$, is:
\[
f
\big(
(x_1,x_2),(y_1,y_2)
\big)
\,=\,
\Big(
\arccosh\,
\big(
1
+
{\textstyle{\frac{(x_1-y_1)^2+(x_2-y_2)^2}{2\,x_2y_2}}}
\big)
\Big)^2,
\]
it comes by complexification
\leqnomode\usetagform{default}
\begin{align}
\label{sqrt-f-c-arcosh}
f^c
\big(
(z_1,z_2),(\overline{z}_1,\overline{z}_2)
\big)
\,=\,
\Big(
\arccosh\,
\big(
1
-
2\,
{\textstyle{\frac{(\Imsmall\,z_1)^2+(\Imsmall\,z_2)^2}{
(\Resmall\,z_2)^2+(\Imsmall\,z_2)^2}}}
\big)
\Big)^2,
\end{align}
with $z_1 = \Re\, z_1 + \isqrt\, \Im\, z_1$ and $z_2 = \Re\, z_2 +
\isqrt\, \Im\, z_2$, provided that
certain inequalities are 
satisfied by $\Im\, z_1$ and $\Im\,z_2$ for this formula to be
meaningful.  Here, the complexification of $\H$ reads as:
\[
\H^c
\,:=\,
\big\{
(z_1,z_2)\in\C^2
\colon\,
\Re\,z_2
>
0
\big\}.
\]

\begin{Lemma}
The domain of definition of $f^c$ in $\H^c$ contains:
\[
\big\{ 
(\Im\, z_1)^2 
<
(\Re\,z_2)^2
\big\}.
\]
\end{Lemma}

\proof
Indeed, the argument $1-2\,Q$ of $\arccosh$ in~(\ref{sqrt-f-c-arcosh}) is real and
$\leqslant 1$. But with $s = \sigma + \isqrt\, t$, for $\cosh\, s$ to
be real $\leqslant 1$, since its imaginary part:
\[
2\,\Im\,
\big(
\cosh\,
s
\big)
\,=\,
2\,\Im\,
\big(
e^{\sigma+it}
+
e^{-\sigma-it}
\big)
\,=\,
\big(e^\sigma-e^{-\sigma}\big)\,
\sin\,t,
\]
vanishes if and only if $t \equiv 0\, \mod\,\pi$, and
since $\cosh\, \sigma > 1$ whenever $\sigma \in \R\backslash\{0\}$,
necessarily $s = \isqrt\, t \in \isqrt\,\R$, hence: 
\[
\arccosh\big(1-2\,Q\big)
\,=:\,
\isqrt\,T
\,\in\,
\isqrt\,\R
\]
for some $T \in \R$, whence:
\[
1-2\,Q
\,=\,
\cosh\,
\big(
\isqrt\,T
\big)
\,=\,
\cos\,T
\eqno
{\scriptstyle{(T\,\in\,\R)}}.
\]

Then $-1 \leqslant \cos\, T \leqslant 1$ forces: 
\[
-\,1
\,\leqslant\,
1
-
2\,
{\textstyle{\frac{(\Imsmall\,z_1)^2+(\Imsmall\,z_2)^2}{
(\Resmall\,z_2)^2+(\Imsmall\,z_2)^2}}}
\,\leqslant\,
1,
\]
the first inequality being equivalent to $(\Im\, z_1)^2 \leqslant
(\Re\, z_2)^2$, while the second holds trivially.
\endproof

For later convenience, let us rewrite the local complex coordinates as
$z_1=u+\isqrt\,v$ and $z_2 = x+ \isqrt\,y$.  Furthermore, let us
restrict our considerations to the subdomain of the above 
domain $\{v^2 \leqslant
x^2\}$ defined by:
\[
0
\,\leqslant\,
1
-
2\,
{\textstyle{\frac{y^2+v^2}{x^2+y^2}}}
\,\leqslant\,
1
\,\,\,\Longleftrightarrow\,\,\,
2y^2+v^2
\,\leqslant\,
x^2,
\]
which guarantees that $\arccosh\, \big( 1 - 2\, \frac{y^2+v^2}{x^2+v^2}
\big)$ is single valued in $[0, \frac{\pi}{2} \big]$.

\begin{center}
\begin{picture}(0,0)%
\includegraphics{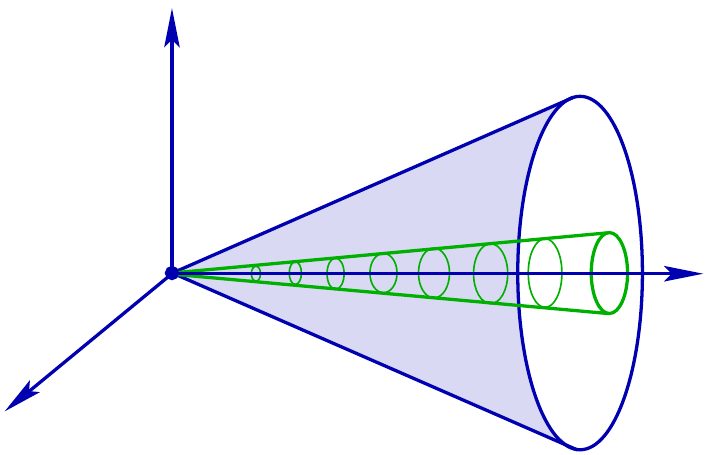}%
\end{picture}%
\setlength{\unitlength}{4144sp}%
\begingroup\makeatletter\ifx\SetFigFont\undefined%
\gdef\SetFigFont#1#2#3#4#5{%
  \reset@font\fontsize{#1}{#2pt}%
  \fontfamily{#3}\fontseries{#4}\fontshape{#5}%
  \selectfont}%
\fi\endgroup%
\begin{picture}(3239,2066)(1464,-2239)
\put(4337,-799){\makebox(0,0)[lb]{\smash{{\SetFigFont{9}{10.8}{\familydefault}{\mddefault}{\updefault}{\color[rgb]{0,0,0}\blue{$\big\{2y^2+v^2\leqslant x^2\big\}$}}%
}}}}
\put(4057,-1188){\makebox(0,0)[lb]{\smash{{\SetFigFont{9}{10.8}{\familydefault}{\mddefault}{\updefault}{\color[rgb]{0,0,0}\green{$\Omega_\varepsilon$}}%
}}}}
\put(1479,-1925){\makebox(0,0)[lb]{\smash{{\SetFigFont{9}{10.8}{\familydefault}{\mddefault}{\updefault}{\color[rgb]{0,0,0}\blue{$y$}}%
}}}}
\put(2287,-284){\makebox(0,0)[lb]{\smash{{\SetFigFont{9}{10.8}{\familydefault}{\mddefault}{\updefault}{\color[rgb]{0,0,0}\blue{$v$}}%
}}}}
\put(4607,-1373){\makebox(0,0)[lb]{\smash{{\SetFigFont{9}{10.8}{\familydefault}{\mddefault}{\updefault}{\color[rgb]{0,0,0}\blue{$x,u$}}%
}}}}
\put(2215,-1547){\makebox(0,0)[lb]{\smash{{\SetFigFont{9}{10.8}{\familydefault}{\mddefault}{\updefault}{\color[rgb]{0,0,0}\blue{$0$}}%
}}}}
\end{picture}%

\end{center}

Drawing $\H = \{x > 0\}$ as a single right half-axis in order to keep
two directions for the $y$- and $v$-axes, this domain $\{2y^2 + v^2 <
x^2\}$ looks like a "security cone" which will contain all subsequent
Grauert tubes $\Omega_\varepsilon$.

Then by the relation:
\[
\arccosh\,(t) 
\,=\,
\isqrt\,\arccos\,t
\eqno
{\scriptstyle{(0\,\leqslant\,t\,\leqslant\,1)}},
\]
we get from~{\thetag{\ref{sqrt-f-c-arcosh}}}
in this subdomain $\{2\,y^2 + v^2 \leqslant x^2\}$ of
$\H^c$:
\[
f^c
\,=\,
-\,
\Big(
\arccos\,
\big(
1
-
2\,
{\textstyle{\frac{y^2+v^2}{x^2+v^2}}}
\big)
\Big)^2,
\]
hence coming back to~(\ref{rho-f-c}):
\leqnomode\usetagform{default}
\begin{align}
\rho
\,=\,
\Big(
\arccos
\big( 
1
-
2\,
{\textstyle{\frac{y^2+v^2}{x^2+y^2}}}
\big)
\Big)^2.
\end{align}

\begin{Lemma}
\label{Omega-varepsilon-pi-2}
For every $0 < \varepsilon < \big( \frac{\pi}{2} \big)^2$, the Grauert
tube around $\H$ in $\H^c$ for the canonical K\"ahler potential
associated with the Poincaré metric on $\H$:
\[
\Omega_\varepsilon
\,:=\,
\big\{
\big(u+\isqrt\,v,\,
x+\isqrt\,y\big)
\in
\H^c
\colon\,
\sqrt{\rho}\,
(u,v,x,y)
<
\sqrt{\varepsilon}
\big\},
\]
has $\mathcal{C}^\omega$ strongly pseudoconvex
boundary $\partial\Omega_\varepsilon = 
\{\rho = \varepsilon\}$ of equation:
\[
2v^2
-
\big(1-\cos\,\sqrt{\varepsilon}\big)\,x^2
+
\big(1+\cos\,\sqrt{\varepsilon}\big)\,y^2
\,=\,
0.
\]
\end{Lemma}

\begin{center}
\includegraphics[scale=0.35]{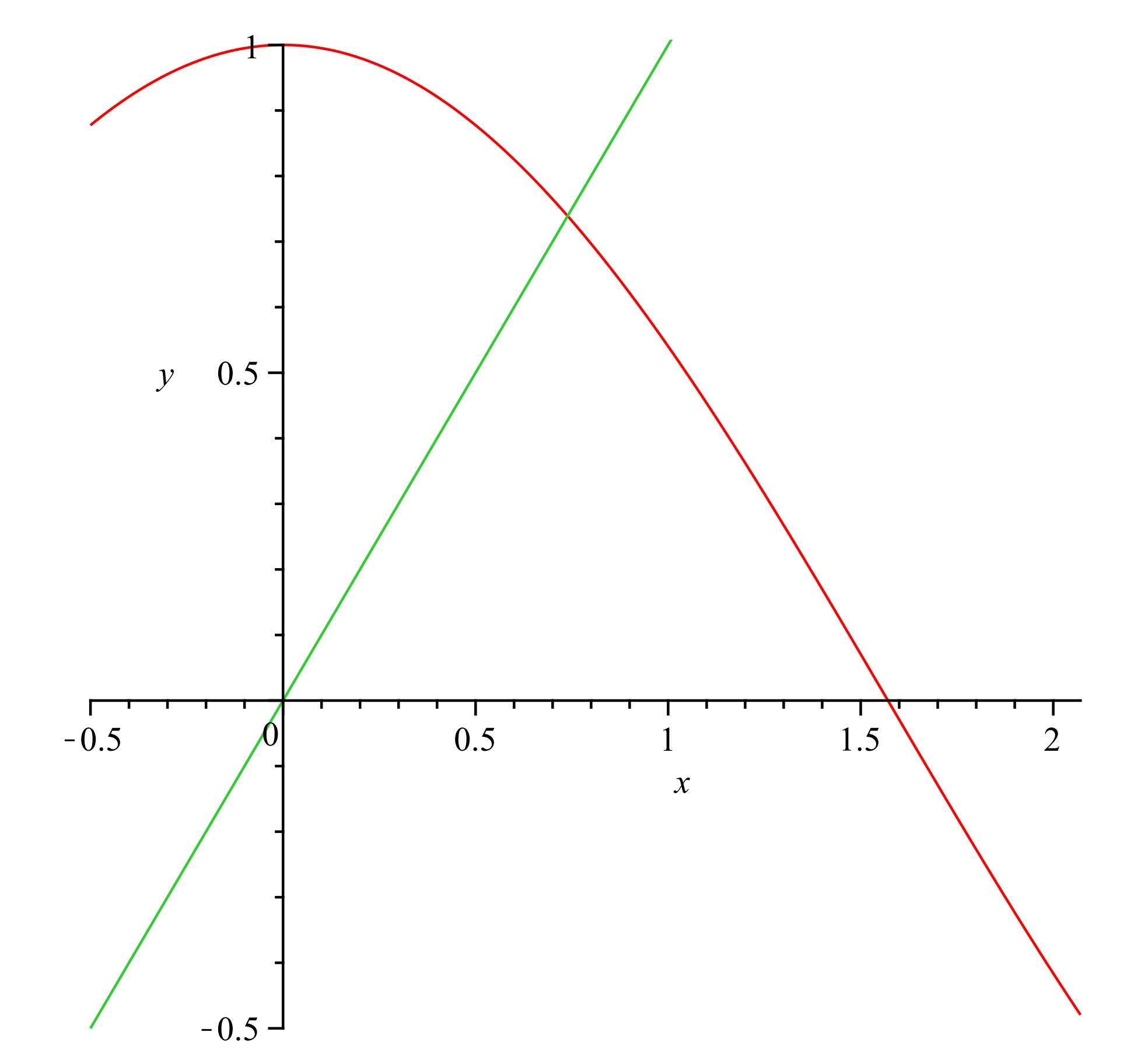}
\ \ \ \ \ \ \ \ \ \ \ \ \ \ \
\includegraphics[scale=0.35]{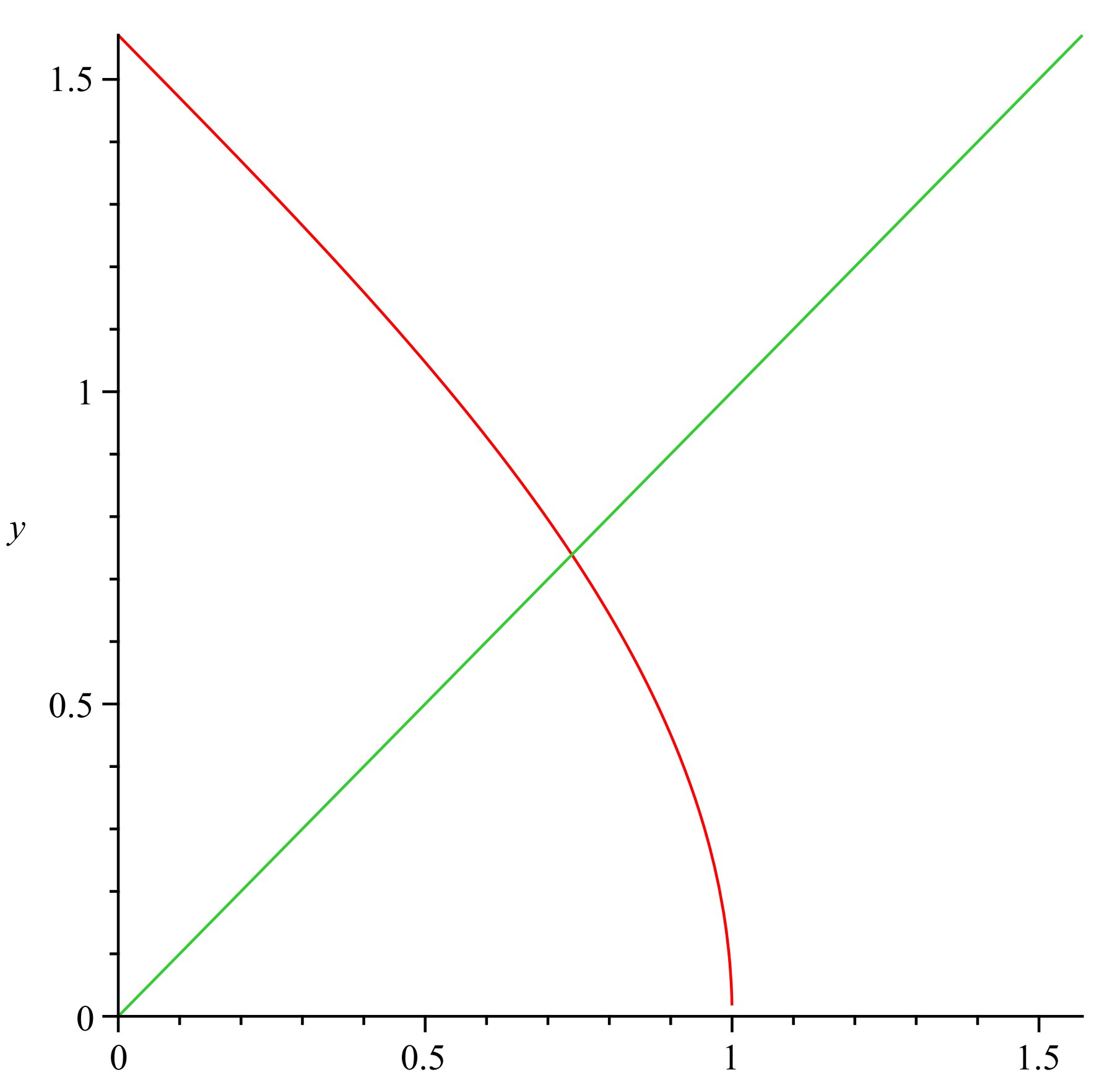}
\end{center}

\proof
Since the function $\arccos$ is a decreasing
$\mathcal{C}^\omega$ diffeomorphism $[0,1) \longrightarrow (0,
\frac{\pi}{2} \big]$, we have:
\[
\aligned
\arccos\,
\big(
1
-
2\,
{\textstyle{\frac{y^2+v^2}{x^2+y^2}}}
\big)
\,<\,
\sqrt{\varepsilon}
&
\,\,\,\Longleftrightarrow\,\,\,
1
-
2\,
{\textstyle{\frac{y^2+v^2}{x^2+y^2}}}
\,>\,
\cos\,\sqrt{\varepsilon}
\\
&
\,\,\,\Longleftrightarrow\,\,\,
\underbrace{
2v^2
-
\big(1-\cos\,\sqrt{\varepsilon}\big)\,x^2
+
\big(1+\cos\,\sqrt{\varepsilon}\big)\,y^2}_{
=:\,r_\varepsilon(u,v,x,y)}
\,<\,
0.
\endaligned
\]

Since $x > 0$, the term $2\,x\,dx$ in the differential
$dr_\varepsilon$ guarantees that $\partial \Omega_\varepsilon = \{
r_\varepsilon = 0\}$ is geometrically smooth at every
point. 

Furthermore, with $w := u + \isqrt\, v$ and $z := x + \isqrt\,
y$, dropping pluriharmonic terms:
\[
r_\varepsilon
\,\equiv\,
w\overline{w}
-
\big(1-\cos\,\sqrt{\varepsilon}\big)\,
{\textstyle{\frac{z\overline{z}}{2}}}
+
\big(1+\cos\,\sqrt{\varepsilon}\big)\,
{\textstyle{\frac{z\overline{z}}{2}}},
\]
we see that $r_\varepsilon$ is strictly plurisubharmonic,
whence $\Omega_\varepsilon = \{ r_\varepsilon < 0\}$ is strongly 
pseudoconvex.
\endproof

In particular, the result holds for thin tubes corresponding
to $0 < \varepsilon \ll \big( \frac{\pi}{2} \big)^2$.

\Section{\bf Calculation of the Complex Cartan Curvature of
$\partial\Omega_\varepsilon \subset \H^c$}
\label{calculation-complex-curvature}

In~\cite{Ebenfelt-Duong-Zaitsev-2016}, the authors proved the
non-existence of CR-umbilical points on the boundaries of Grauert
tubes around flat tori by showing the nonvanishing of a certain
invariant determinant introduced in~\cite{Ebenfelt-Zaitsev-2016},
which vanishes exactly when the Cartan curvatures vanish. In this
paper, we shall use an explicit expression of Cartan curvatures
obtained before by the second named author and Sabzevari
in~\cite{Merker-Sabzevari-2012, Merker-Sabzevari-2014} for locally
graphed hypersufaces, and alternatively a formula
in~\cite{Foo-Merker-Ta-crmath} for hypersurfaces given as zero locus
of implicit functions.

For a $\mathcal{C}^6$-smooth Levi-nondegenerate real 3-dimensional
hypersurface $M \subset \mathbb{C}^2$ represented in complex
coordinates $z = x+\isqrt\, y$, $w = u+\isqrt\, v$ by a local graphing
function:
\[
v
\,=\,
\varphi(x,y,u),
\] 
the Cartan essential curvatures of $M$ are two
real invariants $\mathbf{\Delta_1}$, $\mathbf{\Delta_4}$ expressed
in~\cite[Theorem~1.1]{Merker-Sabzevari-2012}
by following a Tanaka approach, explicitly in terms of
$J_{x,y,u}^6\varphi$, both containing more than 
1,\,500,\,000 terms when expanded.

An equivalent approach~\cite{Merker-Sabzevari-2014} 
closer to Cartan's~\cite{Cartan-1932} can be summarized as
follows. Local generators of $T^{1,0}M$ and $T^{0,1}M$ are:
\[
\mathcal{L}
\,:=\,
\frac{\partial}{\partial z}
-
\frac{\varphi_z}{\isqrt+\varphi_u}\,
\frac{\partial}{\partial u}
\ \ \ \ \ \ \ \ \ \ \ \ \
\text{\rm and}
\ \ \ \ \ \ \ \ \ \ \ \ \
\overline{\mathcal{L}}
\,:=\,
\frac{\partial}{\partial\overline{z}}
-
\frac{\varphi_{\overline{z}}}{-\isqrt+\varphi_u}\,
\frac{\partial}{\partial u},
\]
and their commutator:
\[
\mathcal{T}
\,:=\,
\isqrt\,
\big[
\mathcal{L},
\overline{\mathcal{L}}
\big]
\,=\,
\ell\,
\frac{\partial}{\partial u}
\]
incorporates the {\em real} coefficient, so-called {\sl Levi factor:}
\[
\ell
\,:=\,
2\,\frac{
\varphi_{z\overline{z}}
(1+\varphi_u^2)
-
\isqrt\,\varphi_{\overline{z}}\varphi_{zu}
+
\isqrt\,\varphi_z\varphi_{\overline{z}u}
-
\varphi_{\overline{z}}\varphi_{zu}\varphi_u
-
\varphi_z\varphi_{\overline{z}u}\varphi_u
+
\varphi_z\varphi_{\overline{z}}\varphi_{uu}
}{(1+\varphi_u^2)^2},
\]
which is nowhere vanishing if and only if $M$ is Levi nondegenerate.

Abbreviating the coefficients of $\mathcal{L}$ and 
$\overline{\mathcal{L}}$ as: 
\[
A
\,:=\,
-\,\frac{\varphi_z}{\isqrt+\varphi_u}
\ \ \ \ \ \ \ \ \ \ \ \ \
\text{\rm and}
\ \ \ \ \ \ \ \ \ \ \ \ \
\overline{A}
\,:=\,
-\,\frac{\varphi_{\overline{z}}}{-\isqrt+\varphi_u},
\]
then in terms of the following key function 
(the expansion of which is
$1$ page long):
\[
\overline{P}
\,:=\,
\frac{
\ell_{\overline{z}}
-
\overline{\ell}\,\overline{A}_u
+
\overline{A}\,\ell_u}{
\ell},
\]
the (single) essential Cartan complex invariant expresses in
non-expanded form as:
\begin{align}
\label{invariant-I}
\mathfrak{I}
:=
\frac{1}{6}\,
&\frac{1}{{\sf c}\overline{\sf c}^3}\,
\Big(
-\,2\,
\overline{\mathcal{L}}\big(\mathcal{L}\big(
\overline{\mathcal{L}}\big(\overline{P}\big)\big)\big)
+
\\
&
+3\,\overline{\mathcal{L}}\big(
\overline{\mathcal{L}}\big(
\mathcal{L}\big(\overline{P}\big)\big)\big)
-\,
7\,\overline{P}\,\,
\overline{\mathcal{L}}\big(\mathcal{L}\big(\overline{P}\big)\big) \nonumber
+
\\
&
+
4\,\overline{P}\,
\mathcal{L}\big(\overline{\mathcal{L}}\big(
\overline{P}\big)\big)
-
\mathcal{L}\big(\overline{P}\big)\,
\overline{\mathcal{L}}\big(\overline{P}\big)
+
2\,\overline{P}\,\overline{P}\,\mathcal{L}
\big(\overline{P}\big)
\Big), \nonumber
\end{align}
and a comparison with~\cite{Merker-Sabzevari-2012}
done at the end of~\cite{Merker-Sabzevari-2014} shows 
that it also expresses as:
\[ 
\mathfrak{I}
\,=\,
\frac{4}{
{\sf c}
\overline{\sf c}^3}
\big(
\mathbf{\Delta_1}
+
\isqrt\,\mathbf{\Delta_4}
\big),
\]
where the quantity ${\sf c} \in \C \backslash\{ 0\}$ is a group
parameter of a certain initial $G$-structure, and it has the following
signification.

Suppose there really is a local biholomorphic equivalence $h \colon
\C^2 \longrightarrow \C^2$ which transfers $M$ into $M' := h(M)$, so
that in some appropriate target coordinates $z' = x' + \isqrt\, y'$,
$w' = u' + \isqrt\, v'$, the (localized) image is also graphed as:
\[
v'
\,=\,
\varphi'\big(x',y',u'\big).
\]
Compute similarly $\mathcal{L}'$, $\overline{\mathcal{L}}'$,
$\ell'$, $\overline{P}'$, $\mathfrak{I}'$, but 
extract parts independent of group parameters:
\[
\mathfrak{I}
\,=\,
{\textstyle{\frac{1}{{\sf c}\overline{\sf c}^3}}}\,
\mathfrak{I}_\smallbullet
\ \ \ \ \ \ \ \ \ \ \ \ \
\text{\rm and}
\ \ \ \ \ \ \ \ \ \ \ \ \
\mathfrak{I}'
\,=\,
{\textstyle{\frac{1}{{\sf c'}{\overline{\sf c}'}^3}}}\,
\mathfrak{I}_\smallbullet'.
\]

Because the differential $h_* \colon T\C^2 \longrightarrow T\C^2$
leaves invariant complex tangents, whence $h_*\big( T^{1,0}M\big) =
T^{1,0}M'$, there is a nowhere vanishing function $c' \colon M'
\longrightarrow \C\backslash\{0\}$ such that:
\[
h_\ast\big(\mathcal{L}\big)
\,=\,
c'\,\mathcal{L}'.
\]
At a basic level, it is an easy exercise 
(\cite[p.~44]{Merker-Pocchiola-Sabzevari-2013})
to express the invariancy of the levi factors $\ell$ and $\ell'$ 
through the biholomorphism $h$ as:
\[
\ell
\,=\,
c'\,\overline{c}'\,
\ell',
\]
and at a higher level, a standard feature of Cartan's method of
equivalence then shows that:
\[
\mathfrak{I}_\smallbullet
\,=\,
c'{c'}^3\,
\mathfrak{I}',
\]
which justifies, since $c' \neq 0$ vanishes
nowhere, the invariancy, under changes of
holomorphic coordinates, of the following

\begin{Definition}
A point $p \in M$ at which $\mathfrak{I}(p) = 0$ is called a 
{\sl CR-umbilical point}.
\end{Definition}

In continuation with Lemma~{\ref{Omega-varepsilon-pi-2}} above, 
we are now ready to state and to establish the main proposition.
Inside the complexification of Poincaré's upper half-plane:
\[
\H^c
\,=\,
\big\{
\big(u+\isqrt\,v,
x+\isqrt\,y
\big)
\in
\C^2
\colon\,
x>0
\big\},
\]
consider for every $0 < \varepsilon < \big( \frac{\pi}{2} \big)^2$
the hypersurface:
\[
M_\varepsilon
\,:=\,
\partial\Omega_\varepsilon
\,=\,
\Big\{
\big(
u+\isqrt\,v,
x+\isqrt\,y
\big)
\in
\H^c
\colon\,\,
v^2
-
\big(
{\textstyle{\frac{1-\cossmall\,\sqrt{\varepsilon}}{2}}}
\big)\,
x^2
+
\big(
{\textstyle{\frac{1+\cossmall\,\sqrt{\varepsilon}}{2}}}
\big)\,
y^2
=
0
\Big\}.
\]

\begin{Proposition}
\label{no-CR-umbilics-H-c}
All hypersurfaces $M_\varepsilon \subset \H^c$ 
with $0 < \varepsilon < \big( \frac{\pi}{2} \big)^2$
have {\em no} CR-umbilical point.
\end{Proposition}

\proof
The plain global linear biholomorphism of $\H^c$:
\[
w'
\,:=\,
w,
\ \ \ \ \ \ \ \ \ \ \ \ \
\ \ \ \ \ \ \ \ \ \ \ \ \
z'
\,:=\,
z\,
\sqrt{
{\textstyle{\frac{1+\cossmall\,\sqrt{\varepsilon}}{2}}}
},
\]
transforms $M_\varepsilon$ into:
\[
M_\varepsilon'
\,:=\,
\Big\{
\big(
u'+\isqrt\,v',
x'+\isqrt\,y'
\big)
\in
\H^c
\colon\,\,
{v'}^2
-
{\textstyle{\frac{
1-\cossmall\,\sqrt{\varepsilon}}{
1+\cossmall\,\sqrt{\varepsilon}}}}\,
{x'}^2
+
{y'}^2
=
0
\Big\},
\]
and it is appropriate to set\,\,---\,\,mind the change 
varepsilon $\longmapsto$ epsilon\,\,---:
\[
\epsilon
\,:=\,
\sqrt{
{\textstyle{\frac{
1-\cossmall\,\sqrt{\varepsilon}}{
1+\cossmall\,\sqrt{\varepsilon}}}}
},
\]
so that the equation of $M_\epsilon' := M_\varepsilon'$ 
becomes a bit simpler (dropping the primes):
\[
v^2
-
\epsilon^2\,x^2
+
y^2
\,=\,
0.
\]

Since this fractional map $\varepsilon \longmapsto 
\epsilon(\varepsilon)$ has derivative:
\[
\frac{d}{d\varepsilon}\,
\sqrt{
{\textstyle{\frac{
1-\cossmall\,\sqrt{\varepsilon}}{
1+\cossmall\,\sqrt{\varepsilon}}}}
}
\,=\,
\frac{1}{2\sqrt{\varepsilon}}\,
\frac{\sin(\sqrt{\varepsilon})}{
\sqrt{
{\textstyle{\frac{
1-\cossmall\,\sqrt{\varepsilon}}{
1+\cossmall\,\sqrt{\varepsilon}}}}}\,
\big(1+\cos\,\sqrt{\varepsilon}\big)^2
}
\]
everywhere positive, it is a $\mathcal{C}^\omega$
diffeomorphism $(0, \frac{\pi^2}{4}\big)
\longrightarrow (0,1)$, so that the new $\epsilon$
varies plainly in the open unit real segment:
\[
0
\,<\,
\epsilon
\,<\,
1.
\]

Reminding that $x > 0$, this new equation:
\[
y^2+v^2
\,=\,
\epsilon^2\,x^2,
\]
shows that, a bit similarly as for the flat torus in 
Example~{\ref{flat-torus-example}},
either $v \neq 0$ or $y \neq 0$ at any point. 

Suppose therefore firstly that $v \neq 0$. For the 
$\mathcal{C}^\omega$ graph:
\[
v
\,=\,
\sqrt{\epsilon^2\,x^2-y^2},
\]
a direct calculation of $\mathfrak{I}$ from the formula~(\ref{invariant-I}), by hand or with help of a computer, provides a compact, serendipitous expression:
\[
\mathfrak{I}_\smallbullet
\,=\,
-\,\frac{9}{16}\,
\frac{1-\epsilon^4}{
\big(\epsilon^2x^2-y^2\big)^2}\,
\frac{\big(x+\isqrt\,y\big)^2}{
\big(x-\isqrt\,y\big)^2},
\]
which visibly vanishes nowhere since $x > 0$ whence $(x + \isqrt\,
y)^4 \neq 0$.

Suppose secondly that $y \neq 0$. Since only points with $v = 0$ are
not already examined, assume $v = 0$. For the $\mathcal{C}^\omega$
graph:
\[
y
\,=\,
\sqrt{\epsilon^2\,x^2-v^2}, 
\]
at points with $v = 0$, another direct calculation of the invariant $\mathfrak{I}$ from~(\ref{invariant-I}) also
provides a compact, nowhere vanishing expression:
\[
\mathfrak{I}_\smallbullet
\,=\,
\frac{9}{16}\,
\frac{(1-\epsilon^2)}{
(\epsilon+\isqrt)^2\,\epsilon^4\,x^4},
\]
and this completes the proof of inexistence of CR-umbilical points on $M_\varepsilon$.
\endproof

\proof[Second proof of Proposition~{\ref{no-CR-umbilics-H-c}}]
The formula~(\ref{invariant-I}), explicit as it is, usually gives long
and complicated expression for the combined complex-valued Cartan
invariant $\mathfrak{I}$. This reality is due to the iterated process
of taking roots, derivatives, quotients, {\em etc}. when the graphing
function of the hypersurfaces under consideration is not simple,
including taking roots for example (see, for example, the formulas
given in~\cite{Merker-Sabzevari-2014} and~\cite{ezhov-schmalz}). There
are instances where the hypersurfaces actually have much simpler
representation by mean of implicit functions. An example is the case
of general ellipsoidal hypersurfaces in $\mathbb{C}^2$ considered
in~\cite{Foo-Merker-Ta-crmath}, where a direct calculation from the
formula~(\ref{invariant-I}) for a graphing function of the ellipsoids
gives a very complicated expression for $\mathfrak{I}$, while an
alternative formula (cf. \cite[Corollary 12]{Foo-Merker-Ta-crmath})
applied to simple implicit defining functions of the ellipsoids allows
one to see a whole curve of CR-umbilical points. As the implicit
defining function of $M_\varepsilon$ is also very simple, we shall use
the formulation in~\cite{Foo-Merker-Ta-crmath} to verify the
nonvanishing of the Cartan curvature of $M_\varepsilon$ once again.

Let us recall the necessary formulas
from~\cite{Foo-Merker-Ta-crmath}. For a Levi nondegenerate analytic
hypersurface $M$ in $\mathbb{C}^2$ given by an implicit defining
function:
\[
0
=
F(z,w,\Bar{z},\Bar{w}),
\]
we set
\begin{align*}
    \mathit{L}&:= -F_w \frac{\partial}{\partial z} + F_z \frac{\partial}{\partial w},
    \\
    \mathit{\overline{L}}&:=-F_{\overline{w}} \frac{\partial}{\partial \Bar{z}} + F_{\Bar{z}} \frac{\partial}{\partial \overline{w}},
    \\
    h(F)&:= F_zF_zF_{ww}-2 F_zF_wF_{zw}+ F_wF_wF_{zz},
    \\
    l(F)&:=  F_{\Bar{z}}F_zF_{w\overline{w}} 
            -  F_{\Bar{z}}F_wF_{z\overline{w}}
            -  F_{\overline{w}}F_zF_{\Bar{z}w}
            +  F_{\overline{w}}F_wF_{z\Bar{z}}.
\end{align*}

\begin{Theorem}
(\cite{Foo-Merker-Ta-crmath})
On the domain $\{  F_w \neq 0 \}$, the Cartan invariant $\mathfrak{I}$ of $M$ vanishes exactly on the zero locus of
\begin{align}
\label{invariant-implicit}
  \mathbf{I}_{[w]}
  &
  {\textstyle{
  \, := \,
  12 \,
  \big(F_w\big)^9 \,
  \big( \sum\limits_{i=1}^{7} I_i \big)
  }},
\end{align}
where
\begin{align*}
    I_1
    &=
    {\textstyle{
        \big(
        \frac{l(F)}{F_w^2}
        \big)^3
    }}
    \cdot
    {\textstyle{
        \overline{L}^4
        \big(
        \frac{h(F)}{F_w^3}
        \big)
    }},
    \\
    I_2
    &=
    -6 \, 
    {\textstyle{
        \big( 
        \frac{l(F)}{F_w^2} 
        \big)^2 
    }}
    \cdot
    {\textstyle{
        \mathit{\overline{\mathit{L}}} 
        \big( \frac{l(F)}{F_w^2} 
        \big)
    }}
    \cdot
    {\textstyle{
        \mathit{\overline{\mathit{L}}}^3 
        \big( \frac{h(F)}{F_w^3}
        \big)
    }},
    \\
    I_3
    &=
    -4 \,
    {\textstyle{
        \big( 
        \frac{l(F)}{F_w^2} 
        \big)^2 
    }}
    \cdot
    {\textstyle{
        \mathit{\overline{\mathit{L}}}^2 
        \big( 
        \frac{l(F)}{F_w^2} 
        \big)
    }}
    \cdot
    {\textstyle{
        \mathit{\overline{\mathit{L}}}^2 
        \big( 
        \frac{h(F)}{F_w^3}
        \big)
    }},
    \\
    I_4
    &=
    -
    {\textstyle{
        \big( 
        \frac{l(F)}{F_w^2} 
        \big)^2 
    }}
    \cdot
    {\textstyle{
        \mathit{\overline{\mathit{L}}}^3 
        \big(
        \frac{l(F)}{F_w^2} 
        \big)
    }}
    \cdot
    {\textstyle{
        \mathit{\overline{\mathit{L}}} 
        \big( 
        \frac{h(F)}{F_w^3}
        \big)
    }},
    \\
    I_5
    &=
    15 \,
    {\textstyle{
        \frac{l(F)}{F_w^2} 
    }}
    \cdot
    {\textstyle{
        \Big[\, 
        \mathit{\overline{\mathit{L}}} 
        \big( 
        \frac{l(F)}{F_w^2} 
        \big)
        \Big]^2 
    }}
    \cdot
    {\textstyle{
     \mathit{\overline{\mathit{L}}}^2 
     \big(
     \frac{h(F)}{F_w^3}
     \big)
    }},
    \\
    I_6
    &=
    10 \,
    {\textstyle{
        \frac{l(F)}{F_w^2} 
    }}
    \cdot
    {\textstyle{
        \mathit{\overline{\mathit{L}}} 
        \big( 
        \frac{l(F)}{F_w^2} 
        \big)
    }}
    \cdot
    {\textstyle{
        \mathit{\overline{\mathit{L}}}^2 
        \big( 
        \frac{l(F)}{F_w^2} 
        \big)
    }}
    \cdot
    {\textstyle{  
        \mathit{\overline{\mathit{L}}} 
        \big( 
        \frac{h(F)}{F_w^3}
        \big)
    }},
    \\
    I_7 
    &=
    -15 \,
    {\textstyle{
        \Big[\,
        \mathit{\overline{\mathit{L}}}
        \big( 
        \frac{l(F)}{F_w^2} 
        \big) 
        \Big]^3 
    }}
    \cdot
    {\textstyle{
        \mathit{\overline{\mathit{L}}} 
        \big( 
        \frac{h(F)}{F_w^3}
        \big)
    }}.
\end{align*}
\end{Theorem}

With this formula~(\ref{invariant-implicit}) for checking the
nonvanishing of the Cartan curvature at hand, we now return to our
hypersurface $M_\varepsilon$. We again take advantage of the
elementary biholomorphic transformation as above, and consider the
equivalent model $M'_\epsilon$ whose defining function
writes $v^2-\epsilon^2
x^2+y^2=0,$ with $0 < \epsilon < 1$. 
Switching the notation for coordinates in order
to reach $F_w \neq 0$, namely using instead:
\[
z
\,=\,
u+\isqrt\,v
\ \ \ \ \ \ \ \ \ \ \ \ \
\text{\rm and}
\ \ \ \ \ \ \ \ \ \ \ \ \
w
\,=\,
x+\isqrt\,y,
\]
we can then rewrite:
\[
\aligned
v^2 
- 
\epsilon^2x^2
+
y^2
&
\,=\,
\big(
{\textstyle{\frac{z-\Bar{z}}{2\isqrt}}}
\big)^2 
- 
\epsilon^2
\big(
{\textstyle{\frac{w+\overline{w}}{2}}}
\big)^2 
+ 
\big(
{\textstyle{\frac{w-\overline{w}}{2\isqrt}}}
\big)^2  
\\
&
\,=\,
-\,{\textstyle{\frac{1}{4}}}\,
\big[ 
(z-\Bar{z})^2
+
(1+\epsilon^2)(w^2+\overline{w}^2)
-
2(1-\epsilon^2)\,w\overline{w}
\big]
\\
&
\,=:\,
-\,{\textstyle{\frac{1}{4}}}\,
F(z,w,\Bar{z},\overline{w}),
\endaligned
\]
so that $M_\epsilon' = \{ F = 0\}$, and then as wanted
we have the nowhere vanishing:
\[
F_w 
\,=\,
(1+\epsilon^2)\,2w 
- 
2(1-\epsilon^2)\overline{w}
\,=\,
4\,(\epsilon^2x+\isqrt\,y)
\,\neq\,
0
\]
on $M'_\epsilon$ thanks to our constant assumption $x>0$.  Thus, the
vanishing locus of $\mathbf{I}_{[w]}$ is exactly the set of
CR-umbilical points of $M'_\epsilon$ in this case.

Now, direct calculation from the formula~(\ref{invariant-implicit}),
by hand or preferably on a computer, and keeping in mind that on
$M'_\epsilon$ we always have $v^2 = \epsilon^2x^2-y^2$, gives us:
\[
\aligned
\mathbf{I}_{[w]}
&=
\frac{27}{64}
\epsilon^8
(1-\epsilon^4)
\overline{w}^2
w^6
\\
&=
\frac{27}{64}
\epsilon^8(1-\epsilon^4)
(x-\isqrt y)^2
(x+\isqrt y)^6.
\endaligned
\]

It is then evident that $\mathbf{I}_{[w]}$ is everywhere nonzero on
$M'_\epsilon$ because $x>0$. This completes our second justification
of the inexistence of CR-umbilical points on $M_\varepsilon \cong
M'_\epsilon$.
\proof

\Section{\bf Transfer to Hyperbolic Genus $g \geqslant 2$ Compact 
Surfaces}
\label{transfer-genus-g}

Now, let $S$ be a closed compact oriented $\mathcal{C}^\omega$ surface
of genus $g \geqslant 2$, considered as a Riemann surface.
The Poincar\'e-K\"obe uniformization theorem provides a 
holomorphic covering: 
\[
\tau
\colon\ \ \
\H
\,\longrightarrow\,
S
\,\cong\,
\H
\big/
\pi_1(S).
\]
The Poincaré metric $ds_\H^2 = \lambda\, \big( dx_1^2 + dx_2^2 \big)$
with $\lambda := \frac{1}{x_2^2}$ on $\H$ has constant Gaussian
curvature:
\[
-\,\frac{1}{2\,\lambda}\,
\Big(
\frac{\partial^2}{\partial x_1^2}
+
\frac{\partial^2}{\partial x_2^2}
\Big)
\big(
\log\,\lambda
\big)
\,=\,
-\,1,
\]
and is furthermore kept invariant by all elements
of the group $\Aut\, \H \cong PSL(2,\R)$ 
of holomorphic automorphisms of $\H$:
\[
\big(\Aut\,\H\big)^\ast
\big(ds_\H^2\big)
\,=\,
ds_\H^2,
\]
which acts transitively (and isometrically) on the homogeneous
space $\H$. 

Furthermore, the group of all covering automorphisms of 
$\H \overset{\tau}{\longrightarrow} S$ 
happens to be a discrete {\em subgroup:}
\[
\Aut\,
\Big(
\H
\overset{\tau}{\,\longrightarrow\,}
S
\Big)
\,\subset\,
PSL(2,\R)
\,=\,
\Aut\,\H.
\]
Consequently (and as is well known), 
$ds_\H^2$ descends by push-forward, independently
of preimage points, as a metric on $S$: 
\[
ds_S^2 
\,:=\,
\tau_\ast \big( ds_\H^2\big),
\]
having the same curvature $-1$.

Next, forget the holomorphic structure on $S$, consider now $S$ as
a $\mathcal{C}^\omega$ real surface equipped with this
$\mathcal{C}^\omega$ metric $ds_S^2$, and denote the Bruhat-Whitney
complexification of $S$ by $S^c$. Then
Section~{\ref{Khaler-potential-Grauert-tubes}} gives by
complexification a unique strictly plurisubharmonic $\mathcal{C}^\omega$
 K\"ahler potential $\rho \colon S^c \longrightarrow \R_+$ whose sublevel sets:
\[
\Delta_\varepsilon
\,:=\,
\{
\rho 
<
\varepsilon
\}
\,\subset\,
S^c,
\]
for all small enough $0 < \varepsilon \leqslant \varepsilon_0 \ll 1$,
are strongly pseudoconvex domains bounded by the $\mathcal{C}^\omega$ hypersurfaces:
\[
\partial
\Delta_\varepsilon
\,=\,
\{
\rho
=
\varepsilon
\}.
\]
Here, $\varepsilon_0$ might well be quite small, depending on the
convergence radii of the real-analytic objects that are complexified.

\begin{Lemma}
\label{M-varepsilon-no-CR-umbilic}
Shrinking $\varepsilon_0 > 0$ if necessary, $M_\varepsilon$
has {\em no} CR-umbilical point for all $0 < \varepsilon \leqslant
\varepsilon_0$.
\end{Lemma}

\proof
The uniformizing map, viewed as a $\mathcal{C}^\omega$ map
$\tau \colon \H \longrightarrow S$, also complexifies to become
a holomorphic map:
\[
\H^c
\,\supset\,
V^c
\overset{\tau^c}{\,\longrightarrow\,}
U^c
\,\subset\,
S^c,
\]
where $V^c$ is some open neighborhood of $\H$ in $\H^c$: $\H \subset V^c \subset \H^c$ ,
possibly narrowing much as one reaches $\partial \H = \{x_2 = 0 \}$,
and where $U^c$ is also an open neighborhood of $S$ in $S^c$: $S \subset U^c \subset S^c$.

Since $\tau \colon \H \longrightarrow S$ is a covering map, hence a local
$\mathcal{C}^\omega$ diffeomorphism, each point $p \in S$ has a small
open neighborhood $p \in U_p \subset S$ on which there exist
$\mathcal{C}^\omega$-diffeomorphic inverses of $\tau$, namely maps:
\[
\chi_p
\colon\ \ \
U_p
\overset{\sim}{\,\longrightarrow\,}
\chi_p(U_p)
\,=:\,
V_{\chi_p(p)}
\,\subset\,
\H,
\]
that are uniquely defined as soon as a central point $\chi_p(p) \in
\tau^{-1}(p) \subset \H$ has been chosen in the fiber to fix a level.
Shrinking $U_p$ if necessary, 
the complexification $\chi_p^c$ of $\chi_p(p) $ is also locally biholomorphic
at $p$. 

By compactness of $S \subset S^c$, there 
exists a finite open cover
$U_1^c, \dots, U_\KK^c \subset S^c$ of $S$:
\[
S
\,\subset\,
U_1^c
\cup\cdots\cup
U_\KK^c
\,\subset\,
U^c
\eqno
{\scriptstyle{(\KK\,\geqslant\,1)}},
\]
together with {\em biholomorphic} inverses of the complexification
$\tau^c \colon V^c \longrightarrow U^c$:
\[
\chi_k^c
\colon\ \ \
U_k^c
\overset{\sim}{\,\longrightarrow\,}
\chi_k^c
\big(
U_k^c
\big)
\,=:\,
V_k^c
\,\subset\,
\H^c
\eqno
{\scriptstyle{(1\,\leqslant\,k\,\leqslant\,\KK)}}.
\]

If necessary, shrink $\varepsilon_0 > 0$ so that, 
for all $0 < \varepsilon \leqslant \varepsilon_0$:
\[
\Delta_\varepsilon
\,\subset\,
\Delta_{\varepsilon_0}
\,\Subset\,
U_1^c
\cup\cdots\cup
U_\KK^c.
\]

Now, take any point $q \in \partial \Delta_\varepsilon$. How to
convince oneself that the Cartan CR-curvatures of the strongly pseudoconvex hypersurface $\partial
\Delta_\varepsilon$ is {\em nonzero} at $q$?

This is very simple. For sure, $q \in U_k^c$ for some $1 \leqslant k
\leqslant \KK$. Remind also the tube $\Omega_\varepsilon \subset
\H^c$.  Then because the metric on $S$ is the push-forward of
Poincaré's metric on $\H$, the tubes $\Omega_\varepsilon$ and
$\Delta_\varepsilon$ correspond to each other, namely $\chi_k^c$ sends
$\Delta_\varepsilon \cap U_k^c$ biholomorphically onto
$\Omega_\varepsilon \cap V_k^c$ with:
\[
\chi_k^c\big(q\big)
\,\in\,
\partial\Omega_\varepsilon,
\]
and since the nonvanishing of Cartan CR-curvatures is a biholomorphically
invariant property, Proposition~{\ref{no-CR-umbilics-H-c}}
offers what was wanted.
\endproof

With some basic knowledge on Fuchsian groups, we can also provide a

\proof[Variation on the proof of 
Lemma~{\ref{M-varepsilon-no-CR-umbilic}}]
As already seen, the quotient map:
\[
\tau\colon\ \ \
\mathbb{H} 
\,\longrightarrow\,
S 
\,\cong\, 
\mathbb{H}
\big/
\pi_1(S)
\] 
is locally isometric. Abbreviate: 
\[
G
\,:=\,
\Aut\,
\Big(
\H
\overset{\tau}{\,\longrightarrow\,}
S
\Big)
\,\cong\,
\pi_1(S).
\]

\begin{Definition}
A {\sl fundamental domain} for $S$ is an open subset $D \subset
\mathbb{H}$ whose $G$-translates cover:
\[
\H
\,=\,
\bigcup_{g\in G}\,
g(D), 
\]
being mutually disjoint:
\[
\emptyset
\,=\,
D
\cap
g(D)
\eqno
{\scriptstyle{(\forall\,g\,\in\,G\,\backslash\,\{\Idsmall\})}},
\]
and which has the further property of being
{\sl locally finite} in the sense that 
each compact subset $K \Subset \H$ 
meets only finitely many $G$-images
of $D$.
\end{Definition}

\begin{Theorem}
{\rm (\cite[Chap.~9]{Beardon-1995})}
Relatively compact fundamental domains $D \Subset \H$ having piecewise
$\mathcal{C}^\omega$ boundary consisting of $4g$ geodesic segments
always exist on the universal cover $\tau \colon \H \longrightarrow S$
of any genus $g \geqslant 2$ compact
Riemann surface.\qed
\end{Theorem}

Then in place of a (rough) finite Borel-Lebesgue covering $S \subset
U_1 \cup\cdots\cup U_\KK$ as used in the first proof, we can employ a
geometrically more meaningful 
covering.  For such a fundamental domain $D \subset \H$
of $S$, there is an atlas of $S$ consisting of $4g+1$ open charts:

\smallskip\noindent$\bullet$\, 
$V_0 := D$ itself;

\smallskip\noindent$\bullet$\, 
slightly thickened thin
neighborhoods $V_1, \dots, V_{4g}$ of the $4g$ sides of $D$. 

\smallskip

Further, one can arrange that the restrictions:
\[
\tau
\colon\ \ \
V_i
\,\longrightarrow\,
\tau(V_i)
\,=:\,
U_i
\,\subset\,
S
\eqno
{\scriptstyle{(i\,=\,0,\,1,\,\dots,\,4\,g)}}
\]
are $\mathcal{C}^\omega$ diffeomorphisms. Complexifying
their inverses $\chi_i \colon U_i \overset{\sim}{\longrightarrow}
V_i$ as:
\[
\chi_i^c
\colon\ \ \
U_i^c
\overset{\sim}{\,\longrightarrow\,}
V_i^c
\] 
we can now reason similarly as in the first proof, 
and this concludes.
\endproof

\begin{Remark}
We observe the following interesting facts about the (non)vanishing of
the essential curvatures $\Delta_1$ and $\Delta_4$ on the boundaries
of Grauert tubes of small radii around closed surfaces $S$.

\begin{itemize}

\smallskip\item[{\bf (1)}]\,
If $S$ is the $2$-sphere with the standard round metric, 
both $\Delta_1 $ and
$\Delta_4$ vanish identically.

\smallskip\item[{\bf (2)}]\,
If $S$ is a $2$-dimensional flat torus, we leave as
an exercise to the reader to verify that
$\Delta_1 $ never vanishes, while $\Delta_4$ vanishes
identically.

\smallskip\item[{\bf (3)}]\,
If $S$ is a closed genus $g \geqslant 2$ hyberbolic surface, then both
$\Delta_1 $ and $\Delta_4$ vanish nowhere.

\end{itemize}

\end{Remark}

\Section{\bf Grauert Tubes with Respect to Extrinsic Metrics}
\label{Grauert-extrinsic-metrics}

In Section~{\ref{Khaler-potential-Grauert-tubes}}, Grauert tubes are
constructed with respect to metrics obtained from given {\sl
intrinsic} Riemannian metrics on surfaces. In this section, we look
at constructions of Grauert tubes around surfaces from an {\sl
extrinsic} point of view. More precisely, let us consider a totally
real embedding of a surface $S$ into a complex manifold $X$ of complex
dimension 2. We will identify the surface $S$ with its image under the
embedding, so that $S$ is viewed as a submanifold of $X$. A given
Riemannian metric $d_X$ on $X$ always induces an {\sl extrinsic}
metric on $S$ and the Grauert tubes $\Omega_{\varepsilon}$ around $S$
also can be defined with respect to $d_X$ as $\Omega_{\varepsilon}:=
\{ x \in X: d_X(x,S) < \varepsilon \}$ for small enough positive
$\varepsilon$.

Recall that for a real $n-$dimensional submanifold $M$ of a complex
$n-$dimensional manifold $X$, a point $p$ of $M$ is called a {\sl
complex point} if the tangent vector space of $M$ at $p$ contains at
least one complex line with respect to the complex structure $J \in
\End(TX)$ on the tangent bundle of $X$, that is $T_pM \cap J(T_pM) \neq
\{ 0 \}$. An embedding of $M$ into $X$ is called a {\sl totally real
embedding} if $M$ does not contain any complex point.

It is known that every affine $n$-dimensional totally real vector
subspace $V \subset \C^n$ is affinely holomorphically equivalent to
$\R^n \subset \C^n$. It is also known that every $\mathcal{C}^\omega$
real $n$-dimensional submanifold $M \subset \C^n$ is locally
holomorphically equivalent to $\R^n \subset \C^n$, namely at any point
$p \in M$, there is an open neighborhood $p \in U \subset \C^n$ and a
biholomorphism $h \colon U \overset{\sim}{\longrightarrow} h(U) =: V$
with $h(p) = 0$ such that $h \big( M \cap U \big) = \R^n \cap V$.
Hence an alternative description of maximally real $\mathcal{C}^\omega$
submanifolds $M \subset \C^n$ is as follows.

\begin{Definition}
A real $n$-dimensional $\mathcal{C}^\omega$ 
submanifold $M$ of a complex $n$-dimensional
manifold $X$ is {\sl totally real} if there exists a 
family indexed by $\alpha \in A$ of biholomorphisms: 
\[
\varphi_\alpha 
\colon\ \ \
U_\alpha 
\overset{\sim}{\,\longrightarrow\,}
\varphi_\alpha\big(U_\alpha\big)
\,=:\,
V_\alpha
\,\subset\,
X
\]
with $U_\alpha \subset \C^n$ open, with $V_\alpha \subset X$ open,
with $X = \medcup_\alpha V_\alpha$, such that:

\smallskip\noindent$\bullet$\,
if $\varphi_\alpha(0) \not\in M$, then $\varphi_\alpha\big(U_\alpha\big)
\cap M = \emptyset$;

\smallskip\noindent$\bullet$\,
if $\varphi_\alpha(0) \in M$, then the restriction:
\[
\varphi_\alpha
\big\vert_{\R^n\cap U_\alpha}
\colon\ \ \ \ \
\R^n\cap U_\alpha
\overset{\sim}{\,\longrightarrow\,}
M\cap V_\alpha,
\]
is a $\mathcal{C}^\omega$ real diffeomorphism.
\end{Definition}

\begin{Example}
By looking at the standard complex atlas of the complex projective
space $\mathbb{CP}^n$, it is clear that $\mathbb{RP}^n$ is totally
real in $\mathbb{CP}^n$. On $\mathbb{RP}^n$, there is a canonical
round metric induced from the round metric on its double cover
$\mathbb{S}^n$. The Guillemin-Stenzel metric associated to this round
metric on $\mathbb{CP}^n$ is nothing but the Fubini-Study metric
on $\mathbb{CP}^n$. The
complexified manifold $(\mathbb{S}^n)^c$ is a double cover of
$\mathbb{CP}^n$, which is a real $2n-$dimensional submanifold in the
$\mathbb{S}^1-$fibration $\mathbb{S}^{2n+1}$ of $\mathbb{CP}^n$.
\end{Example}
 
\begin{Example}
Of particular interest for us here is the fact that a product of two
totally real submanifolds is also totally real, which is evident from
either definition.
\end{Example}
 
\begin{Example}
\label{Example-flat-products}
Let us look at Example 3.3 once again, this time from an extrinsic
point of view. Consider a 2-dimensional real vector subspace $V$ of
$\C^2$ which passes through the origin, with coordinates $(z, w) \in
\C^2$. The intersections of $V$ with the $z$-axis and $w$-axis are two
real line. Therefore $V$ can be written in exactly one of the
following three forms.

\smallskip\noindent{\footnotesize\sf Case~1:}
$V = \big\{ y = \alpha\,x,\, v = \beta\, u \big\}$, 
where $\alpha, \beta$ are real. The Grauert
tube $\Omega_\varepsilon(V)$ of radius $\varepsilon$ around $V$ with
respect to the standard distance in $\C^2$ is given by:
\[
\Big\{
(x + \isqrt\, y, u + \isqrt\, v)
    \in 
    \C^2: 
        \frac{
            (\alpha x - y)^2
                }{
                    (\alpha^2+1)^2
                    } 
        + \frac{
            (\beta u - v)^2
                }{
                (\beta^2+1)^2
                } 
 < \varepsilon^2 
 \Big\}.
 \]

In order to obtain a compact hypersurface, we take the quotient
$\widetilde{\Omega}_\varepsilon(V)$ of $\Omega_\varepsilon(V)$ by the
translations by $2\pi$ on each real coordinates of $V$. Then
$\widetilde{\Omega}_\varepsilon(V)$ can be embedded into $\C^2$ as:
\[
\Big\{ 
(z,w)\in\C^2
\colon\,
\Big( 
\log\big\vert
e^{\frac{\isqrt\,z}{1+\isqrt\, \alpha}}
\big\vert
\Big)^2 
+ 
\Big( 
\log\,
\big\vert
e^{ \frac{\isqrt\, w}{ 1 + \isqrt\, \beta}}
\big\vert
\Big)^2 
< 
\varepsilon^2 
\Big\}.
\]
Any point on the boundary of $\widetilde{\Omega}_\varepsilon(V)$
admits the same local defining function as its preimage on the
boundary of $\Omega_\varepsilon(V)$. Solving the local defining
function for the variable $v$ gives the graph:
\[
    v = 
    \beta u 
    - ({\beta^2+1}) 
        \sqrt{ 
            \varepsilon^2 -  
            \frac{(\alpha x - y)^2
            }{
                (\alpha^2+1)^2
                }
        }.
\]

A direct calculation of the Cartan invariant
using the formula~({\ref{invariant-I}})
provides:
\[ 
\mathfrak{J}_{\smallbullet} = 
\frac{
-9\,(\alpha + \isqrt)^9\,
(-\alpha + \isqrt)^{11}\,
(\beta + \isqrt)^{16}\,
(-\beta + \isqrt)^{16}\,
\varepsilon^8
}{  
( -\alpha x + y + \varepsilon + a^2 \varepsilon)^8\,
( \alpha x - y + \varepsilon + a^2 \varepsilon)^8},
\]
and this result is nowhere vanishing. 
So the boundary of $\widetilde{\Omega}_\varepsilon(V)$ 
also does not contain any CR-umbilical point.

\smallskip\noindent{\footnotesize\sf Case~2:}
$V = \big\{ x = 0,\, v = \beta u \big\}$, where $\beta$ is again real.
The Grauert tube of radius $\varepsilon$ around $V$ with respect to
the standard distance in $\C^2$ is now given by:
\[ 
\Omega_\varepsilon(V)
\,=\, 
\Big\{ 
(x + \isqrt\, y, u + \isqrt\, v)
\in 
\C^2: 
x^2 + \frac{(\beta u - v)^2}{(\beta^2+1)^2} 
< 
\varepsilon^2 
\Big\}.
\]

A point on the boundary of $\widetilde{\Omega}_\varepsilon(V)$ or of
${\Omega}_\varepsilon(V)$ admits the local graphing function:
\[
v
\,=\, 
\beta u 
-
(\beta^2+1)
\sqrt{\varepsilon^2- x^2},
\]
of which the (relative) Cartan curvature can be computed from the
formula~(\ref{invariant-I}) to be:
\[
\mathfrak{J}_{\smallbullet} 
\,=\,
\frac{9\,(\beta^2+1)^{16}\,\varepsilon^8
}{
(x^2-\varepsilon^2)^8}.
\]
Thus, the (relative) invariant $\mathfrak{J_{\smallbullet}}$ is also nowhere vanishing on the boundary. 

Note that $\widetilde{\Omega}_\varepsilon(V)$ can be embedded into 
$\C^2$ as:
\[
\Big\{
(z,w)
\in
\C^2
\colon\,
\big(\log\,\vert e^{z}\vert\big)^2 
+ 
\Big(
\log\,\big\vert e^{\frac{\isqrt\,w}{1+\isqrt\,\beta}}\big\vert
\Big)^2 
< 
\varepsilon^2 
\Big\}.
\]
 
\smallskip\noindent{\footnotesize\sf Case~3:}
$V = \big\{ x = 0 = u \big\}$. A point on
$\widetilde{\Omega}_\varepsilon(V)$ 
which can be embedded into $\C^2$ as:
\[
\big\{
(z,w)\in\C^2
\colon\,
\big(\log\,\vert e^{z}\vert\big)^2
 + 
\big(\log\,\vert e^{w}\vert\big)^2 
< 
\varepsilon^2 
\big\},
\]
now admits the local defining function 
\[ 
x^2+u^2 
\,=\,
\varepsilon^2.
\]

In this case, we do not obtain a local graphing function of the form
$v = \phi(x,y,u)$, but a simple calculation using the alternative
formula~(\ref{invariant-implicit}) for the implicit defining function
$F(z,w, \Bar{z}, \overline{w}) =
(\frac{z + \Bar{z}}{2})^2 + (\frac{w + \overline{w}}{2})^2 - 
\varepsilon^2$ 
shows that the relative invariant
$\mathfrak{J_{\smallbullet}}$ is proportional to:
\[ 
\frac{27\,(x^2+u^2)^4}{64} 
\,=\,
\frac{27\,\varepsilon^{8}}{64}.
\]
So it is evident that the boundary of
$\widetilde{\Omega}_\varepsilon(V)$ also does not contain any
CR-umbilical point.
\end{Example}

For two given Riemannian manifolds $(X,d_X),(Y,d_Y)$, the distance
$d_{X \times Y}$ with respect to the product metric on $X
\times Y$ is:
\begin{align}
d^2_{X \times Y}\big((x_1,y_1),\,(x_2,y_2)\big)
\,=\,
d^2_{X}(x_1,x_2)
+
d^2_{Y}(y_1,y_2),
\end{align}
assuming that $X,Y$ are uniquely geodesic, {\em i.e.} there exists a
unique geodesic between any two points. 

Our next examples of Grauert tubes in $\C \times \C$ will be
constructed with respect to products of two extrinsic metrics on $\C
\supset \R$. For the two possible component metrics on $\C$, we will
consider the three standard ones: flat, elliptic and hyperbolic.

\smallskip\noindent$\bullet$\,
{\footnotesize\sf Flat metric on $\C$.}
Denote by $d_{\Flatsmall}$
the flat Pythagorean 
metric on $\C \ni x + \isqrt\, y$. Consider the
totally real line $V_{\Flatsmall} = \{ y=0 \}$ in
$U_{\Flatsmall} = \C$. The flat distance from
any point $z \in U_{\Flatsmall}$ to $V_{\Flatsmall}$ is:
\begin{align}
d_{\Flatsmall}\big(z,\,V_{\Flatsmall}\big) 
\,=\,
\big\vert\Im(z)\big\vert
\,=\,
\vert y\vert.
\end{align}

\smallskip\noindent$\bullet$\,
{\footnotesize\sf Elliptic metric on $\C\P^1$.}
For the elliptic metric $d_{\Ellsmall}$, we look at the local chart
$U_0=\{ [1:z]: z\in \C \}$ of $\mathbb{CP}^1.$ Since
$(\mathbb{CP}^1,d_{\Ellsmall})$ is not uniquely geodesic, we consider
a small neighborhood $U_{\Ellsmall}=\{ [1:z]: |z| < \delta \}$ of
$[1:0]$ in $U_0$, which is uniquely geodesic for small positive
$\delta$ thanks to the fact that the injective radius of
$(\mathbb{CP}^1,d_{\Ellsmall})$ is positive. Then $V_{\Ellsmall}=\{
[1:\Re(z)]: \Re(z) < \delta \}$ is totally real in $U_{\Ellsmall}.$
 
\begin{Lemma}
The elliptic distance from any point $(x,y) \approx [1:(x+\isqrt\,
y)] $ of $ U_{\Ellsmall}$ to $V_{\Ellsmall}$ is given by:
\begin{align}
\,\,\,\,\,\,\,\,\,\,\,
d_{\Ellsmall}\big((x,y),\,V_{\Ellsmall}\big) 
\,=\,
\arccos
\Big(
\frac{\sqrt{1+x^2}}{\sqrt{1+x^2+y^2}}
\Big).
\end{align}
\end{Lemma}
 
\begin{proof}
A point $[1:(x + \isqrt\, y)]$ of $\mathbb{CP}^1$ corresponds to the
point $\big( 
\frac{1}{ \sqrt{1+x^2 + y^2}}, \frac{x}{ \sqrt{1+x^2 + y^2}},
\frac{y}{ \sqrt{1+x^2 + y^2}} \big)$ 
of $\mathbb{S}^2$ embedded in $\R^3$,
and a point $[1: \alpha]$ of $V_{\Ellsmall}$ corresponds to 
$\big( \frac{1}{
\sqrt{1+\alpha^2}}, \frac{\alpha}{ \sqrt{1+\alpha^2}},0 \big)$.
 
Now $d_{\Ellsmall} \big( (x,y), V_{\Ellsmall} \big)$ 
is exactly the spherical
distance between $ P = \big( \frac{1}{ \sqrt{1+x^2 + y^2}}, \frac{x}{
\sqrt{1+x^2 + y^2}}, \frac{y}{ \sqrt{1+x^2 + y^2}} \big)$ and the arc
$\{ Q_{\alpha}=(\frac{1}{ \sqrt{1+\alpha^2}}, \frac{\alpha}{
\sqrt{1+\alpha^2}},0) : \alpha>0 \}$, that is:
\[
\cos\, 
d_{\Ellsmall}\big((x,y),\,V_{\Ellsmall}\big)
\,=\,
\underset{\alpha>0}{\max}\,
\frac{\langle P,Q_{\alpha}\rangle}{|P|\,|Q_{\alpha}|}.
\]
Using the Cauchy-Schwartz inequality, we have:
\begin{align*}
\frac{\langle P,Q_{\alpha}\rangle
}{|P|\,|Q_{\alpha}|}
&
= 
\frac{1+\alpha x}{\sqrt{1+\alpha^2}\,\sqrt{1+x^2+y^2}} 
\\
&
\leqslant
\frac{\sqrt{1+\alpha^2}\,\sqrt{1+x^2}}{
\sqrt{1+\alpha^2}\,\sqrt{1+x^2+y^2}} 
\\
&
=
\frac{\sqrt{1+x^2}}{\sqrt{1+x^2+y^2}},               
\end{align*}
where the maximum is attained at $\alpha = x$.  
\end{proof}

\smallskip\noindent$\bullet$\,
{\footnotesize\sf Hyperbolic metric on $\H$.}
For the hyperbolic metric $d_{\Hypsmall}$, we may
consider a small open
neighborhood $U$ of 0 in the Poincar\'e disc, and the totally real
interval $U \cap \{ \Im(z)=0 \}$ in $U$,
but it is more convenient to work
with the corresponding domain $U_{\Hypsmall}=\{ z= x+ \isqrt\, y\}$ of
$U$ on the upper-half plane model, which is an open neighborhood of
$\isqrt\,$. The corresponding totally real interval in $U_{\Hypsmall}$
is $V_{\Hypsmall} = U_{\Hypsmall} \cap \{ \Re(z)=0 \}$.

\begin{Lemma}
The hyperbolic distance from any point $(x,y) \approx z = x+\isqrt\, 
y$ in $U_{\Hypsmall}$ to $V_{\Hypsmall}$ is given by:
\begin{align}
d_{\Hypsmall}
\big((x,y),\,V_{\Hypsmall}\big) 
\,=\,
\arccosh\,
\Big(
\frac{\sqrt{x^2+y^2}}{y}
\Big).
\end{align}
\end{Lemma}

\proof
Recall that for a hyperbolic triangle on the upper-half plane with
angles $A$, $B$, $C$ and opposite sides of lengths $a$, $b$, $c$, 
the rule of sine reads:
\[ 
\frac{\sin A}{\sinh\,a} 
= 
\frac{\sin B}{\sinh\,b} 
=
\frac{\sin C}{\sinh\,c}.
\]

Thus, given the angle $A$ and the side $a$, the side $b$ is of maximal
length when $B = \frac{\pi}{2}$ because the function $\sinh$ is
monotone and because:
\[ 
\sinh\,b 
= 
\sin B 
\;
\frac{\sinh\,a}{\sin\,A} 
\leqslant 
\frac{\sinh\,a}{\sin\,A}. 
\]

It follows that to find the hyperbolic distance from a given point
$z = x+\isqrt\, y $ to the line $V_{\Hypsmall}$, 
we look at the geodesic
line passing through $z$ and orthogonal to $V_{\Hypsmall}$, which is
the half-circle on the upper-half plane model with centre at 0 and of
radius $|z| = \sqrt{x^2+y^2}$. This geodesic line intersects
$V_{\Hypsmall}$ at the point $\big( 0,\sqrt{x^2+y^2}
\big) \approx 0+\isqrt\,
\sqrt{x^2+y^2}$. Thus, we have:
\begin{align}
d_{\Hypsmall}\big( 
(x,y),\,V_{\Hypsmall} 
\big) 
&
=
d_{\Hypsmall}
\big( 
(x,y),\,(0,\sqrt{x^2+y^2})
\big)
\notag
\\
&
=
\arccosh
\Big( 1+ \frac{(x-0)^2+(y-\sqrt{x^2+y^2})^2}{2y\sqrt{x^2+y^2}}
\Big)
\notag  
\\
&
=
\arccosh
\Big(
\frac{\sqrt{x^2+y^2}}{y}
\Big).
\qedhere
\end{align}
\endproof

We are now in position to give some non-trivial examples of Grauert
tubes with respect to extrinsic metrics.
 
\begin{Proposition}
The Grauert tubes of radius $\varepsilon$ with respect to the product
metric $d_1 \times d_2$ around the totally real submanifold $V_1
\times V_2$ in $U_1 \times U_2$ admit local defining
functions:
\[
\aligned
\rho(x,y,u,v) 
\,:=\,
&\,
\big[
d_1\big((x,y),V_1\big)
\big]^2
+ 
\big[
d_2\big((u,v),V_2\big) 
\big]^2
\\
\,<\,
&\,
\varepsilon^2,
\endaligned
\]
where $(U_i,V_i,d_i)$ for $i=1,2$ is one of the three models:
\[
\big(U_{\Flatsmall},\,V_{\Flatsmall},d_{\Flatsmall}\big),
\ \ \ \ \ \ \ \ \
\big(U_{\Ellsmall},V_{\Ellsmall},d_{\Ellsmall}\big),
\ \ \ \ \ \ \ \ \
\big(U_{\Hypsmall},V_{\Hypsmall},d_{\Hypsmall}\big).
\eqno\qed
\]
\end{Proposition}
 
In particular, we obtain six examples of Grauert tubes with respect to
the corresponding extrinsic product metrics.
 
\begin{Remark} 
Notice here that our examples are of local nature, and not
compact. When both $d_1$ and $d_2$ are flat metrics, one recovers the
local graphing function of the flat torus as in Example 3.3, since:
\[
\big[
d_{\Flatsmall}\big((x,y),V_{\Flatsmall}\big)
\big]^2
+ 
\big[
d_{\Flatsmall}\big((u,v),V_{\Flatsmall}\big)
\big]^2
\,=\,
y^2 + v^2. 
\]
However, the remaining five examples are
very different from those obtained from intrinsic metrics in Example
3.1, Example 3.3 and Lemma 4.4. Thus, the Grauert tubes around the
same totally real manifolds with respect to intrinsic and extrinsic
metrics look very different.
\end{Remark}

\begin{Lemma}
In terms of:
\[
H
\,:=\,
\sqrt{\varepsilon^2
-
\Big[\arccosh \frac{\sqrt{x^2+y^2}}{y}\Big]^2}
\]
and of:
\[
E
\,:=\,
\sqrt{\varepsilon^2 
- 
\Big[\arccos\,\frac{\sqrt{1+x^2}}{\sqrt{1+x^2+y^2}}\Big]^2},
\]
the local defining functions for the boundaries of the Grauert tubes of
radius $\varepsilon$ with respect to the product metrics are given
by Table~1.
\end{Lemma}

\begin{table}
\caption{}\label{eqtable}
\renewcommand\arraystretch{2}
\noindent\[
\begin{array}{|c|c|}
\hline
{\text{Product metrics}}&{\text{Defining functions}}\\
\hline
{d_{\Flatsmall} \oplus d_{\Flatsmall} }&{v =\sqrt{\varepsilon^2 - y^2}}\\
\hline
{d_{\Ellsmall} \oplus d_{\Flatsmall} }&{v =(\varepsilon^2 - \arcsin \, \frac{y}{\sqrt{1+x^2+y^2}}})^{1/2}\\
\hline
{d_{\Hypsmall} \oplus d_{\Flatsmall} }&{v =(\varepsilon^2 - \arcsinh \, \frac{x}{y}})^{1/2}\\
\hline
{d_{\Hypsmall} \oplus d_{\Hypsmall} }&{v = \frac{u}{\sinhsmall H}}\\
\hline
{d_{\Ellsmall} \oplus d_{\Hypsmall} }&{v = \frac{u}{\sinhsmall E}}\\
\hline
{d_{\Ellsmall} \oplus d_{\Ellsmall} }&{v = \frac{ 1 + \sqrt{1 - 4(1+u^2)(\sinsmall E)^2}}{2 \; \sinsmall E}}\\
\hline
\end{array}
\]
\end{table}

\proof
We only treat the case of the product between the
hyperbolic and flat metrics, in which the local graphing function is
given by:
\begin{align}
\rho(x,y,u,v)
\,=\,
\bigg[ 
\arccosh\,
\Big(
\frac{\sqrt{x^2+y^2}}{y}
\Big) 
\bigg]^2
+
v^2
\,<\,
\varepsilon^2,
\end{align}
while the calculations for the other cases can be done in a similar
way.

The defining function for the boundary of the Grauert tube is obtained
by solving the equation $\rho = \varepsilon^2 $ for the variable $v$
as follows:
\begin{align*}
\rho = \varepsilon^2 
&\implies 
\arccosh
\Big(
\frac{\sqrt{x^2+y^2}}{y}
\Big) = \sqrt{ \varepsilon^2 - v^2} 
\\
&\implies
\frac{\sqrt{x^2+y^2}}{y} 
= 
\cosh\,\big(\sqrt{\varepsilon^2-v^2}\big) 
\\
&\implies
1+\frac{x^2}{y^2} 
= 
\big[\cosh\,\big(\sqrt{\varepsilon^2-v^2}\big)\big]^2
=
\big[\sinh\,\big(\sqrt{\varepsilon^2-v^2}\big)\big]^2 +1 
\\
&\implies
\frac{x}{y} 
= 
\sinh\big(\sqrt{\varepsilon^2-v^2}\big). 
\end{align*}

So, the defining function belongs to the rigid case with the graph:
\[
v 
\,=\,
\sqrt{\varepsilon^2-\arcsinh\,
\big(
{\textstyle{\frac{x}{y}}}
\big)}.
\qedhere
\]
\endproof

Unfortunately, except for the case of $d_{\Flatsmall} \oplus
d_{\Flatsmall}$, the expressions of the Cartan invariant obtained by
calculations with either formula~(\ref{invariant-I})
or~(\ref{invariant-implicit}), though explicit, are overwhelmingly
complicated, and so do not allows us to see the CR-umbilical locii.


\bibliographystyle{amsplain}


\end{document}